\documentclass{amsart}
\usepackage{amsthm,amsmath,amssymb,amscd,graphics,enumerate, stmaryrd,xspace,verbatim, epic, eepic}
\usepackage[all]{xypic}
\SelectTips{cm}{}

{
   \newtheorem{theorem}[subsubsection]{Theorem}
        
   \newtheorem{lemma}[subsubsection]{Lemma}
   
      \newtheorem*{fact}{Fact}

   \newtheorem{conjecture}[subsubsection]{Conjecture}
   
}
{\theoremstyle{definition}

   \newtheorem{example}[subsubsection]{Example}
   \newtheorem{definition}[subsubsection]{Definition}
   \newtheorem{remark}[subsubsection]{Remark}
   
}
\newcommand{\CC}{{\mathbb{C}}}
\newcommand{\QQ}{{\mathbb{Q}}}

\newcommand{\PP}{{\mathbb{P}}}
\newcommand{\ZZ}{{\mathbb{Z}}}

\newcommand{\GG}{{\mathbb{G}}}
\newcommand{\LL}{{\mathbb{L}}}

\renewcommand{\AA}{{\mathbb{A}}}

\newcommand{\bR}{{\mathbf{R}}}

\newcommand{\bmu}{{\boldsymbol{\mu}}}

\newcommand{\cB}{{\mathcal B}}
\newcommand{\cC}{{\mathcal C}}
\newcommand{\cE}{{\mathcal E}}

\newcommand{\cG}{{\mathcal G}}

\newcommand{\cI}{{\mathcal I}}
\newcommand{\cK}{{\mathcal K}}

\newcommand{\cM}{{\mathcal M}}

\newcommand{\cO}{{\mathcal O}}

\newcommand{\cX}{{\mathcal X}}
\newcommand{\cY}{{\mathcal Y}}

\newcommand{\fC}{{\mathfrak C}}
\newcommand{\fD}{{\mathfrak D}}
\newcommand{\fM}{{\mathfrak M}}

\newcommand{\tw}{{\operatorname{tw}}}

\newcommand{\vir}{{\text{vir}}}
\def\<{\langle}
\def\>{\rangle}

\newcommand{\Isom}{\operatorname{Isom}}

\newcommand{\Sing}{\operatorname{Sing}}

\newcommand{\Hom}{{\operatorname{Hom}}}

\newcommand{\Aut}{{\operatorname{Aut}}}

\newcommand{\lrar}{\longrightarrow}
\newcommand{\dar}{\downarrow}

\newcommand{\ocI}{\overline{{\mathcal I}}}

\newcommand{\ocM}{\overline{{\mathcal M}}}

\newcommand{\double}{\genfrac..{0pt}1
{\raise -1pt\hbox{$\scriptstyle\longrightarrow$}}{\raise 3pt\hbox
{$\scriptstyle\longrightarrow$}}} 

\newcommand{\q}{\bold{q}}

\def\tototi{\mathbin{\mathop{\otimes}\limits^{\raise-1pt\hbox
{$\scriptscriptstyle {\rm L}$}}}}

\def\indlim{\mathop{\vrule width0pt height7pt depth
4pt\smash{\lim\limits_{\raise 1pt\hbox to 14.5pt
{\rightarrowfill}}}}}
\def\projlim{\mathop{\vrule width0pt height7pt depth
4pt\smash{\lim\limits_{\raise 1pt\hbox to 14.5pt
{\leftarrowfill}}}}}

\newcommand\displaceamount{3pt}

\newcommand{\doubledown}{\ar@<\displaceamount>[d]\ar@<-\displaceamount>[d]}

\newcommand{\doubleup}{\ar@<\displaceamount>[u]\ar@<-\displaceamount>[u]}

\newcommand{\doubleright}{\ar@<\displaceamount>[r]\ar@<-\displaceamount>[r]}

\newcommand{\iX}[1][]{{\cI_{{#1}}(\cX)}}
\newcommand{\riX}[1][]{{\overline{\cI}_{{#1}}(\cX)}}

\newcommand{\thickslash}{\mathbin{\!\!\pmb{\fatslash}}}

\newcommand{\lecture}[1]{\section{#1}}

\begin{document}

\title{Lectures on Gromov--Witten invariants of orbifolds}

\author[D. Abramovich]{Dan Abramovich}
\address{Department of Mathematics, Box 1917, Brown University,
Providence, RI, 02912, U.S.A} 
\email{abrmovic@math.brown.edu}


\maketitle

\setcounter{tocdepth}{1}
\tableofcontents

\section*{Introduction}

\subsection{What this is} This text came out of my CIME
minicourse at Cetraro, June 6-11, 2005. I kept the text relatively
close to what actually happened in the course. In particular, because
of last minutes changes in schedule, the lectures on usual
Gromov--Witten theory started after I gave two lectures, so I decided to
give a sort of introduction to non-orbifold Gromov--Witten
theory, including an exposition of Kontsevich's formula for rational plane curves.  From here the gradient of difficulty is relatively high, but I still hope different readers of rather spread-out backgrounds will get something out of it. I gave few computational examples at the
end, partly because of lack of time. An additional lecture was given by Jim Bryan on his work on the crepant resolution conjecture with
Graber and with Pandharipande, with what I find very exciting computations, and I make some comments on this in the last lecture. In a way, this is the original reason for the existence of  orbifold Gromov--Witten theory. 

\subsection{Introspection}
One of the organizers' not--so--secret reasons for inviting me to give
these lectures was to push me and my collaborators to finish the paper
\cite{AGV}. The organizers were only partially successful: the paper,
already years overdue, was only being circulated in rough form upon
demand at the time of these lectures. Hopefully it will be made fully
available by the time these notes are published. Whether they meant it
or not, one of the outcomes of this  intention of the organizers is
that these lectures are centered completely around the paper
\cite{AGV}. I am not sure I did wisely by focussing so much on our
work and not bringing in approaches and beautiful applications so many
others have contributed, but this is the outcome and I hope it serves
well enough. 

\subsection{Where does all this come from?} Gromov--Witten
invariants of orbifolds first appeared, in response to string
theory - in particular Zaslow's pioneering \cite{Zaslow} - in the inspiring work of W. Chen and Y. Ruan (see, e.g.,
\cite{Chen-Ruan}). The special case of orbifold cohomology of finite 
quotient orbifolds was treated in a paper of
Fantechi-G\"ottsche \cite{Fantechi-Gottsche}, and the special case of
that for symmetric product orbifolds was also discovered by Uribe
\cite{Uribe}. In the algebraic setting, the underlying construction of
moduli spaces was undertaken in \cite{AVfibered}, \cite{AV}, and the
algebraic analogue of the work of Chen-Ruan was worked out in
\cite{AGV0} and the forthcoming \cite{AGV}. The case of  finite
quotient orbifolds was developed in  work of
Jarvis--Kauffmann--Kimura \cite{Jarvis-Kaufmann-Kimura}. 
Among the many applications I note the work of Cadman \cite{Cadman},
which has a theoretical importance later in these lectures. 

In an amusing twist of events, Gromov--Witten invariants of orbifolds
were destined to be introduced by Kontsevich but failed to do so - I
tell the story in the appendix. 

\subsection{Acknowledgements}
First I'd like to thank my collaborators Tom Graber and Angelo
Vistoli, whose joint work is presented here, hopefully adequately.
I thank the organizers of the summer school - Kai Behrend, Barbara
Fantechi and Marco Manetti -  for inviting me to give
these lectures. Thanks to  Barbara Fantechi and Damiano Fulghesu whose notes taken at  Cetraro were of great value for the preparation of this text. Thanks are due to Maxim Kontsevich and Lev Borisov who
gave permission to include their correspondence.
Finally, it is a pleasure to acknowledge that this particular project
was inspired by the aforementioned work of W. Chen and Ruan. 

\renewcommand{\thesection}{Lecture \arabic{section}}
\renewcommand{\thesubsection}{\arabic{section}.\arabic{subsection}}

\lecture{Gromov--Witten theory}

\subsection{Kontsevich's formula}

Before talking of orbifolds, let us step back to the story of
Gromov--Witten invariants. Of course these first came to be famous due
to their role in mirror symmetry. But this failed to be excite me, a
narrow-minded algebraic geometer such as I am, until Kontsevich
\cite{Kontsevich} gave his formula for the number of rational plane
curves.  

This is a piece of magic which I will not resist describing.
\subsubsection{Setup}
Fix an integer $d>0$. Fix points $p_1, \ldots p_{3d-1}$ in general
position in the plane. Look at the following number: 
\begin{definition}
$$N_d = {\mbox{\huge$\#$}} 
   \left\{
     \begin{array}{l} C \subset \PP^2 \mbox{   a rational curve,}\\
 \deg C = d, \mbox{ and} \\p_1, \ldots p_{3d-1} \in C
 \end{array}
  \right\}.$$
\end{definition} 
\begin{remark} One sees that $3d-1$ is the right number of points using an elementary dimension count: a degree $d$ map of $\PP^1$ to the plane is parametrized by three forms of degree $d$ (with $3(d+1)$  parameters). Rescaling the forms (1 parameter) and automorphisms of $\PP^1$ (3 parameters) should be crossed out, giving $3(d+1) - 1 - 3 \ = \ 3d-1$.
\end{remark}

\subsubsection{Statement}
\begin{theorem}[Kontsevich] For $d>1$ we have
$$N_d=\sum_{\begin{array}{c}d=d_1+d_2\\d_1, d_2>0\end{array}}  N_{d_1}N_{d_2}\left(d_1^2d_2^2{{3d-4}\choose {3d_1-2}}-d_1^3d_2{{3d-4}\choose{3d_1-1}} \right).$$
\end{theorem}

\begin{remark} The first few numbers are
$$N_1 = 1, \ N_2 =1,\ N_3 = 12,\ N_4 = 620,\ N_5 = 87304.$$
The first two are elementary, the third is classical, but $N_4$ and $N_5$ are nontrivial.
\end{remark}
\begin{remark} The first nontrivial analogous number in $\PP^3$ is the
  number of lined meeting four other lines in general position (the
  answer is 2, which is the beginning of Schubert calculus). 
\end{remark}

\subsection{Set-up for a streamlined proof}
\subsubsection{$\ocM_{0,4}$}
We need one elementary moduli space: the compactified space of ordered four-tuples of points on a line, which we describe in the following unorthodox manner : 
$$\ocM_{0,4} = \overline{\left\{ p_1,p_2,q,r \in L\  \Big|\  
\begin{array}{l} L \simeq \PP^1\\  p_1,p_2,q,r  \mbox{  distinct }
\end{array}
\right\}}$$
The open set $\cM_{0,4}$ indicated in the braces is isomorphic to $\PP^1 \smallsetminus \{0, 1, \infty\}$, the coordinate corresponding to the cross ratio $$CR(p_1, p_2, q,r) = \frac{p_1 - p_2}{p_1 - r}\frac{q-r}{q-p_2}.$$ 
The three points in the compactification, denoted 
\begin{align*} 0& = (p_1,p_2\ |\ q,r), \\
                            1& = (p_1,q\ |\ p_2,r), \mbox{ and}\\ 
                            \infty& = (p_1,r\ |\ p_2,q), \end{align*}
%
\setlength{\unitlength}{0.00083333in}
\begingroup\makeatletter\ifx\SetFigFont\undefined%
\gdef\SetFigFont#1#2#3#4#5{%
  \reset@font\fontsize{#1}{#2pt}%
  \fontfamily{#3}\fontseries{#4}\fontshape{#5}%
  \selectfont}%
\fi\endgroup%
{\renewcommand{\dashlinestretch}{30}
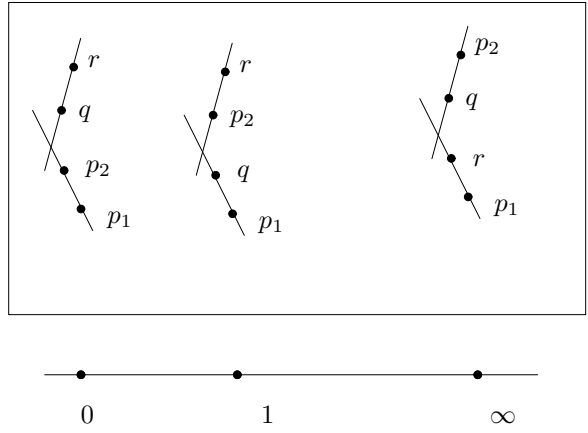
\begin{center}\begin{picture}(3624,2710)(0,-10)
\path(237,358)(3312,358)
\put(462,58){\makebox(0,0)[lb]{{0}}}
\put(1587,58){\makebox(0,0)[lb]{{1}}}
\put(3012,58){\makebox(0,0)[lb]{{$\infty$}}}
\put(1437,358){\blacken\ellipse{50}{50}}
\put(1437,358){\blacken\ellipse{50}{50}}
\put(2937,358){\blacken\ellipse{50}{50}}
\put(2937,358){\blacken\ellipse{50}{50}}
\put(462,358){\blacken\ellipse{50}{50}}
\put(462,358){\blacken\ellipse{50}{50}}
\put(1302,1603){\blacken\ellipse{50}{50}}
\put(1407,1363){\blacken\ellipse{50}{50}}
\put(1287,1978){\blacken\ellipse{50}{50}}
\put(1362,2248){\blacken\ellipse{50}{50}}
\put(2772,1708){\blacken\ellipse{50}{50}}
\put(2877,1468){\blacken\ellipse{50}{50}}
\put(2757,2083){\blacken\ellipse{50}{50}}
\put(2832,2353){\blacken\ellipse{50}{50}}
\put(357,1633){\blacken\ellipse{50}{50}}
\put(462,1393){\blacken\ellipse{50}{50}}
\put(342,2008){\blacken\ellipse{50}{50}}
\put(417,2278){\blacken\ellipse{50}{50}}
\path(12,2683)(3612,2683)(3612,733)
	(12,733)(12,2683)
\path(1407,2428)(1182,1603)
\path(1107,1978)(1482,1228)
\path(2877,2533)(2652,1708)
\path(2577,2083)(2952,1333)
\path(462,2458)(237,1633)
\path(162,2008)(537,1258)
\put(1572,1243){\makebox(0,0)[lb]{{$p_1$}}}
\put(1392,1903){\makebox(0,0)[lb]{{$p_2$}}}
\put(1452,2248){\makebox(0,0)[lb]{{$r$}}}
\put(3042,1348){\makebox(0,0)[lb]{{$p_1$}}}
\put(627,1273){\makebox(0,0)[lb]{{$p_1$}}}
\put(492,1588){\makebox(0,0)[lb]{{$p_2$}}}
\put(447,1933){\makebox(0,0)[lb]{{$q$}}}
\put(507,2278){\makebox(0,0)[lb]{{$r$}}}
\put(1437,1558){\makebox(0,0)[lb]{{$q$}}}
\put(2862,2008){\makebox(0,0)[lb]{{$q$}}}
\put(2922,2353){\makebox(0,0)[lb]{{$p_2$}}}
\put(2907,1663){\makebox(0,0)[lb]{{$r$}}}
\end{picture}\end{center}}
describing the three different ways to split the four points in two pairs and position them on a nodal curve with two rational components.

\subsubsection{A one-parameter family}
We now look at our points $p_1,\ldots p_{3d-1}$ in the plane. 

We pass two lines $\ell_1,\ell_2$ with general slope through the last point $p_{3d-1}$ and consider the following family of rational plane curves in $C \to B$ parametrized by a curve $B$:

\begin{itemize}
\item Each curve $C_b$ contains $p_1,\ldots p_{3d-2}$ (but not necessarily $p_{3d-1}$).
\item One point $q \in C_b \cap \ell_1$ is marked.
\item One point $r \in C_b \cap \ell_2$ is marked.
\end{itemize}

In fact, we have a family of rational curves $C \to B$ parametrized by $B$, most of them smooth, but finitely many have a single node, and a morphism $f: C \to \PP^2$ immersing the fibers in the plane.

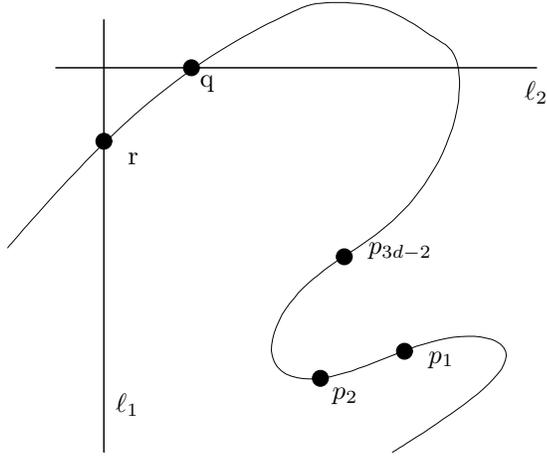
\setlength{\unitlength}{0.00083333in}
\begingroup\makeatletter\ifx\SetFigFont\undefined%
\gdef\SetFigFont#1#2#3#4#5{%
  \reset@font\fontsize{#1}{#2pt}%
  \fontfamily{#3}\fontseries{#4}\fontshape{#5}%
  \selectfont}%
\fi\endgroup%
{\renewcommand{\dashlinestretch}{30}
\begin{picture}(3836,2846)(0,-10)
\put(612,1952){\blacken\ellipse{100}{100}}
\put(1157,2412){\blacken\ellipse{100}{100}}
\put(2112,1232){\blacken\ellipse{100}{100}}
\put(1962,472){\blacken\ellipse{100}{100}}
\put(2487,642){\blacken\ellipse{100}{100}}
\path(312,2412)(3312,2412)
\path(312,2412)(3312,2412)
\path(612,2712)(612,12)
\path(612,2712)(612,12)
\path(12,1287)(13,1288)(14,1290)
	(17,1293)(21,1298)(28,1306)
	(36,1316)(47,1328)(60,1344)
	(76,1362)(94,1383)(115,1407)
	(138,1433)(163,1462)(191,1493)
	(220,1526)(251,1560)(283,1596)
	(317,1633)(352,1671)(388,1710)
	(425,1749)(463,1789)(502,1829)
	(541,1870)(582,1910)(624,1951)
	(666,1992)(710,2033)(755,2074)
	(802,2116)(850,2157)(899,2199)
	(950,2241)(1003,2283)(1058,2325)
	(1113,2367)(1171,2408)(1229,2448)
	(1287,2487)(1355,2529)(1420,2568)
	(1482,2604)(1539,2635)(1591,2663)
	(1638,2687)(1680,2708)(1716,2726)
	(1748,2741)(1775,2754)(1798,2765)
	(1817,2774)(1834,2781)(1848,2787)
	(1860,2792)(1871,2796)(1881,2800)
	(1891,2802)(1901,2805)(1912,2807)
	(1924,2810)(1939,2812)(1955,2813)
	(1974,2815)(1997,2817)(2022,2818)
	(2051,2819)(2084,2819)(2121,2818)
	(2160,2816)(2202,2813)(2247,2807)
	(2292,2799)(2337,2787)(2378,2773)
	(2416,2757)(2452,2739)(2486,2721)
	(2516,2703)(2545,2686)(2570,2669)
	(2594,2653)(2615,2639)(2635,2626)
	(2652,2614)(2668,2604)(2683,2594)
	(2697,2586)(2710,2578)(2722,2570)
	(2733,2563)(2743,2556)(2753,2548)
	(2763,2539)(2772,2530)(2781,2519)
	(2790,2506)(2798,2491)(2805,2474)
	(2812,2454)(2818,2430)(2824,2403)
	(2828,2372)(2830,2337)(2830,2297)
	(2828,2253)(2824,2205)(2815,2152)
	(2803,2096)(2787,2037)(2769,1985)
	(2748,1932)(2725,1879)(2699,1827)
	(2671,1777)(2641,1729)(2611,1682)
	(2579,1639)(2547,1597)(2514,1559)
	(2480,1522)(2447,1488)(2413,1456)
	(2379,1427)(2345,1399)(2311,1372)
	(2277,1347)(2244,1323)(2210,1300)
	(2177,1278)(2143,1256)(2110,1234)
	(2077,1212)(2044,1190)(2012,1167)
	(1980,1143)(1948,1118)(1917,1093)
	(1886,1066)(1856,1037)(1827,1007)
	(1800,975)(1773,942)(1749,908)
	(1726,872)(1706,834)(1689,796)
	(1675,758)(1665,719)(1659,682)
	(1658,646)(1662,612)(1674,577)
	(1693,547)(1718,522)(1747,502)
	(1781,487)(1817,477)(1857,471)
	(1897,469)(1939,471)(1982,476)
	(2026,483)(2070,494)(2114,506)
	(2159,520)(2203,536)(2248,552)
	(2292,569)(2337,587)(2382,605)
	(2427,623)(2472,640)(2517,657)
	(2562,673)(2608,687)(2654,701)
	(2700,712)(2746,722)(2791,729)
	(2836,734)(2880,736)(2922,735)
	(2962,732)(3000,725)(3034,715)
	(3063,702)(3087,687)(3107,665)
	(3119,641)(3122,615)(3118,588)
	(3108,560)(3092,531)(3070,501)
	(3044,471)(3014,440)(2981,409)
	(2944,377)(2905,346)(2864,314)
	(2821,282)(2778,251)(2733,220)
	(2689,190)(2647,162)(2606,135)
	(2567,110)(2533,88)(2502,68)
	(2476,51)(2454,38)(2438,28)
	(2426,20)(2418,16)(2414,13)(2412,12)
\put(1212,2262){\makebox(0,0)[lb]{{q}}}
\put(762,1812){\makebox(0,0)[lb]{{r}}}
\put(2262,1212){\makebox(0,0)[lb]{{$p_{3d-2}$}}}
\put(2037,312){\makebox(0,0)[lb]{{$p_2$}}}
\put(2637,537){\makebox(0,0)[lb]{{$p_1$}}}
\put(687,237){\makebox(0,0)[lb]{{$\ell_1$}}}
\put(3237,2187){\makebox(0,0)[lb]{{$\ell_2$}}}
\end{picture}
}

\subsubsection{The geometric equation}
We have a cross-ratio map 
\begin{align*}B &\stackrel{\lambda}{\lrar} \ocM_{0,4}\\
    C  & \mapsto CR(p_1,p_2,q,r)
    \end{align*}

Since points on $\PP^1$ are homologically equivalent we get
$$\deg_B \lambda^{-1}(p_1,p_2|q,r) = \deg_B \lambda^{-1}(p_1,q|p_2,r).$$

\subsubsection{The right hand side}
Now, each curve counted in $\deg_B \lambda^{-1}(p_1,q|p_2,r)$ is of the following form: 
\begin{itemize}
\item It has two components $C_1,C_2$ of respective degrees $d_1,d_2$ satisfying $d_1+d_2 = d$.
\item We have $p_1\in C_1$ as well as $3d_1-2$ other points among the $3d-4$ points $p_3,\ldots p_{3d-2}$.
\item We have $p_2\in C_2$ as well as the remaining $3d_2-2$  points from $p_3,\ldots p_{3d-2}$.
\item We select one point $z \in C_1\cap C_2$. where the two abstract curves are attached. 
\item We mark one point $q \in C_1 \cap \ell_1$ and one point $r\in C_2\cap \ell_2$. 
\end{itemize}

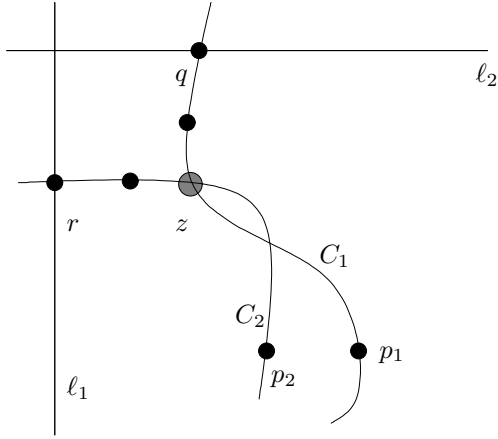
\setlength{\unitlength}{0.00083333in}
\begingroup\makeatletter\ifx\SetFigFont\undefined%
\gdef\SetFigFont#1#2#3#4#5{%
  \reset@font\fontsize{#1}{#2pt}%
  \fontfamily{#3}\fontseries{#4}\fontshape{#5}%
  \selectfont}%
\fi\endgroup%
{\renewcommand{\dashlinestretch}{30}
\begin{picture}(3536,2739)(0,-10)
\put(1212,2412){\blacken\ellipse{100}{100}}
\put(312,1587){\blacken\ellipse{100}{100}}
\put(2207,537){\blacken\ellipse{100}{100}}
\put(1632,537){\blacken\ellipse{100}{100}}
\put(782,1597){\blacken\ellipse{100}{100}}
\put(1137,1962){\blacken\ellipse{100}{100}}
\texture{44555555 55aaaaaa aa555555 55aaaaaa aa555555 55aaaaaa aa555555 55aaaaaa 
	aa555555 55aaaaaa aa555555 55aaaaaa aa555555 55aaaaaa aa555555 55aaaaaa 
	aa555555 55aaaaaa aa555555 55aaaaaa aa555555 55aaaaaa aa555555 55aaaaaa 
	aa555555 55aaaaaa aa555555 55aaaaaa aa555555 55aaaaaa aa555555 55aaaaaa }
\put(1157,1577){\shade\ellipse{150}{150}}
\path(12,2412)(3012,2412)
\path(12,2412)(3012,2412)
\path(312,2712)(312,12)
\path(312,2712)(312,12)
\path(1287,2712)(1286,2710)(1285,2705)
	(1283,2696)(1280,2682)(1275,2663)
	(1269,2638)(1262,2606)(1254,2569)
	(1244,2527)(1234,2480)(1224,2429)
	(1213,2376)(1202,2321)(1191,2266)
	(1181,2211)(1171,2156)(1162,2103)
	(1154,2052)(1148,2003)(1142,1956)
	(1138,1912)(1135,1870)(1133,1831)
	(1132,1794)(1133,1760)(1135,1728)
	(1138,1698)(1143,1670)(1149,1643)
	(1157,1618)(1165,1594)(1176,1571)
	(1187,1549)(1202,1526)(1218,1503)
	(1236,1482)(1255,1461)(1276,1440)
	(1299,1420)(1323,1401)(1349,1382)
	(1376,1363)(1404,1345)(1434,1327)
	(1464,1309)(1496,1292)(1527,1275)
	(1560,1258)(1593,1242)(1625,1225)
	(1658,1209)(1690,1193)(1722,1177)
	(1753,1161)(1784,1145)(1813,1129)
	(1841,1113)(1869,1097)(1895,1081)
	(1920,1064)(1943,1047)(1966,1030)
	(1987,1012)(2010,990)(2032,968)
	(2052,944)(2071,919)(2089,893)
	(2106,865)(2121,837)(2136,808)
	(2149,778)(2161,747)(2172,716)
	(2181,685)(2190,653)(2197,622)
	(2203,590)(2208,560)(2212,530)
	(2215,501)(2218,473)(2219,446)
	(2219,420)(2219,396)(2218,373)
	(2216,351)(2214,331)(2212,312)
	(2208,287)(2202,263)(2196,243)
	(2188,224)(2179,206)(2168,190)
	(2156,174)(2142,159)(2126,146)
	(2110,132)(2093,120)(2077,110)
	(2063,101)(2052,95)(2044,90)
	(2039,88)(2037,87)
\path(87,1587)(89,1587)(94,1587)
	(102,1588)(115,1588)(134,1589)
	(159,1591)(191,1592)(228,1594)
	(272,1595)(321,1597)(375,1599)
	(432,1601)(493,1602)(556,1603)
	(620,1604)(684,1605)(748,1605)
	(811,1605)(872,1604)(931,1603)
	(988,1601)(1042,1598)(1093,1595)
	(1141,1591)(1186,1586)(1228,1581)
	(1268,1574)(1305,1567)(1339,1559)
	(1371,1551)(1400,1541)(1427,1530)
	(1453,1519)(1476,1506)(1498,1492)
	(1518,1478)(1537,1462)(1556,1443)
	(1574,1423)(1590,1401)(1605,1377)
	(1617,1352)(1628,1324)(1638,1294)
	(1646,1262)(1652,1227)(1658,1189)
	(1661,1148)(1664,1104)(1665,1058)
	(1665,1008)(1663,956)(1661,901)
	(1657,843)(1653,784)(1648,724)
	(1642,664)(1636,604)(1629,547)
	(1622,491)(1616,440)(1610,394)
	(1604,353)(1599,319)(1595,291)
	(1592,270)(1590,254)(1588,245)
	(1587,239)(1587,237)
\put(387,237){\makebox(0,0)[lb]{{$\ell_1$}}}
\put(2937,2187){\makebox(0,0)[lb]{{$\ell_2$}}}
\put(1962,1062){\makebox(0,0)[lb]{{$C_1$}}}
\put(1437,687){\makebox(0,0)[lb]{{$C_2$}}}
\put(2337,462){\makebox(0,0)[lb]{{$p_1$}}}
\put(1662,312){\makebox(0,0)[lb]{{$p_2$}}}
\put(1062,1287){\makebox(0,0)[lb]{{$z$}}}
\put(387,1287){\makebox(0,0)[lb]{{$r$}}}
\put(1062,2187){\makebox(0,0)[lb]{{$q$}}}
\end{picture}
}

For every choice of splitting $d_1+d_2 = d$ we have ${3d-4}\choose {3d_1-2}$ ways to choose the set of $3d_1-2$ points on $C_1$ from the $3d-4$ points $p_3,\ldots p_{3d-2}$. We have $N_{d_1}$ choices for the curve $C_1$ and $N_{d_2}$ choices for $C_2$. We have $d_1\cdot d_2$ choices for $z$, $d_1$ choices for $q$ and $d_2$ for $r$. This gives the term

$$\deg_B \lambda^{-1}(p_1,q|p_2,r)\ =  \ \sum_{\begin{array}{c}d=d_1+d_2\\d_1, d_2>0\end{array}} \  {{3d-4}\choose {3d_1-2}}\ \  \cdot\ \   N_{d_1}\ N_{d_2}\ \  \cdot\ \   d_1 d_2\ \  \cdot\ \   d_1\cdot d_2.$$

A simple computation in deformation theory shows that each of these curves actually occurs in a fiber of the family $C \to B$, and it occurs exactly once with multiplicity 1.

\subsubsection{The left hand side}

Curves counted in $\deg_B \lambda^{-1}(p_1,p_2|q,r)$ come in two flavors: there are {\em irreducible} curves passing through $q=r= \ell_1\cap\ell_2$ This is precisely $N_d$.

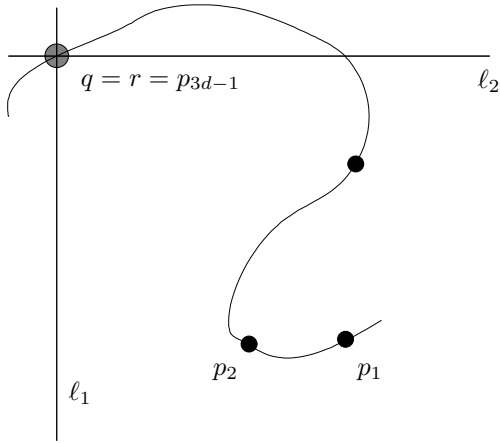
\setlength{\unitlength}{0.00083333in}
\begingroup\makeatletter\ifx\SetFigFont\undefined%
\gdef\SetFigFont#1#2#3#4#5{%
  \reset@font\fontsize{#1}{#2pt}%
  \fontfamily{#3}\fontseries{#4}\fontshape{#5}%
  \selectfont}%
\fi\endgroup%
{\renewcommand{\dashlinestretch}{30}
\begin{picture}(3542,2760)(0,-10)
\texture{44555555 55aaaaaa aa555555 55aaaaaa aa555555 55aaaaaa aa555555 55aaaaaa 
	aa555555 55aaaaaa aa555555 55aaaaaa aa555555 55aaaaaa aa555555 55aaaaaa 
	aa555555 55aaaaaa aa555555 55aaaaaa aa555555 55aaaaaa aa555555 55aaaaaa 
	aa555555 55aaaaaa aa555555 55aaaaaa aa555555 55aaaaaa aa555555 55aaaaaa }
\put(318,2412){\shade\ellipse{150}{150}}
\put(2183,1737){\blacken\ellipse{100}{100}}
\put(1518,612){\blacken\ellipse{100}{100}}
\put(2118,642){\blacken\ellipse{100}{100}}
\path(18,2412)(3018,2412)
\path(18,2412)(3018,2412)
\path(318,2712)(318,12)
\path(318,2712)(318,12)
\path(18,2037)(18,2038)(17,2039)
	(17,2042)(16,2046)(15,2052)
	(14,2060)(13,2070)(12,2081)
	(12,2093)(12,2107)(13,2122)
	(16,2138)(20,2154)(25,2171)
	(33,2189)(43,2208)(56,2227)
	(72,2246)(91,2267)(116,2289)
	(145,2312)(179,2335)(219,2360)
	(266,2386)(318,2412)(361,2432)
	(405,2451)(450,2469)(493,2487)
	(534,2504)(573,2520)(609,2535)
	(642,2549)(672,2563)(699,2575)
	(724,2587)(746,2599)(767,2610)
	(786,2620)(803,2630)(820,2640)
	(837,2650)(853,2659)(871,2668)
	(889,2676)(909,2685)(931,2693)
	(955,2700)(982,2708)(1012,2714)
	(1046,2720)(1084,2725)(1126,2729)
	(1171,2732)(1220,2733)(1273,2731)
	(1328,2728)(1385,2721)(1443,2712)
	(1497,2700)(1548,2686)(1598,2671)
	(1644,2655)(1688,2639)(1728,2622)
	(1766,2607)(1800,2592)(1831,2578)
	(1859,2565)(1885,2553)(1908,2542)
	(1930,2532)(1950,2522)(1968,2513)
	(1986,2504)(2002,2496)(2018,2487)
	(2033,2478)(2048,2468)(2063,2458)
	(2078,2447)(2094,2434)(2109,2419)
	(2125,2403)(2141,2384)(2158,2363)
	(2174,2339)(2191,2312)(2207,2282)
	(2223,2249)(2237,2212)(2249,2172)
	(2259,2129)(2266,2084)(2268,2037)
	(2266,1989)(2259,1941)(2248,1893)
	(2234,1848)(2217,1806)(2198,1766)
	(2177,1729)(2155,1696)(2132,1665)
	(2108,1638)(2084,1614)(2059,1591)
	(2034,1572)(2008,1554)(1983,1537)
	(1957,1522)(1931,1507)(1905,1493)
	(1880,1479)(1853,1465)(1827,1450)
	(1800,1434)(1774,1416)(1747,1397)
	(1719,1376)(1692,1352)(1664,1326)
	(1637,1298)(1609,1266)(1582,1232)
	(1556,1195)(1530,1156)(1505,1115)
	(1482,1072)(1461,1029)(1443,987)
	(1426,941)(1413,899)(1404,860)
	(1397,826)(1393,795)(1391,768)
	(1391,745)(1392,726)(1395,710)
	(1398,697)(1402,686)(1407,677)
	(1412,670)(1418,664)(1424,660)
	(1431,656)(1437,652)(1445,649)
	(1453,645)(1462,641)(1471,636)
	(1481,630)(1493,624)(1505,616)
	(1519,608)(1535,598)(1552,588)
	(1571,577)(1592,566)(1616,555)
	(1641,545)(1668,537)(1704,530)
	(1741,527)(1778,528)(1816,531)
	(1853,537)(1889,545)(1925,555)
	(1961,567)(1996,580)(2031,594)
	(2065,610)(2099,626)(2132,642)
	(2164,659)(2195,676)(2223,692)
	(2250,707)(2273,720)(2293,732)
	(2310,742)(2323,750)(2332,755)
	(2338,759)(2342,761)(2343,762)
\put(393,237){\makebox(0,0)[lb]{{$\ell_1$}}}
\put(2943,2187){\makebox(0,0)[lb]{{$\ell_2$}}}
\put(468,2187){\makebox(0,0)[lb]{{$q=r=p_{3d-1}$}}}
\put(1293,387){\makebox(0,0)[lb]{{$p_2$}}}
\put(2193,387){\makebox(0,0)[lb]{{$p_1$}}}
\end{picture}
}

Now, each {\em reducible} curve counted in $\deg_B \lambda^{-1}(p_1,p_2|q,r)$ is of the following form: 
\begin{itemize}
\item It has two components $C_1,C_2$ of respective degrees $d_1,d_2$ satisfying $d_1+d_2 = d$.
\item We have  $3d_1-1$  points among the $3d-4$ points $p_3,\ldots p_{3d-2}$ are on $C_1$.
\item We have $p_1,p_2\in C_2$ as well as the remaining $3d_2-2$  points from $p_3,\ldots p_{3d-2}$.
\item We select one point $z \in C_1\cap C_2$. where the two abstract curves are attached. 
\item We mark one point $q \in C_1 \cap \ell_1$ and one point $r\in C_1\cap \ell_2$. 
\end{itemize}

\setlength{\unitlength}{0.00083333in}
\begingroup\makeatletter\ifx\SetFigFont\undefined%
\gdef\SetFigFont#1#2#3#4#5{%
  \reset@font\fontsize{#1}{#2pt}%
  \fontfamily{#3}\fontseries{#4}\fontshape{#5}%
  \selectfont}%
\fi\endgroup%
{\renewcommand{\dashlinestretch}{30}
\begin{picture}(3536,2739)(0,-10)
\put(312,1962){\blacken\ellipse{100}{100}}
\put(837,2412){\blacken\ellipse{100}{100}}
\put(1312,2112){\blacken\ellipse{100}{100}}
\put(1737,687){\blacken\ellipse{100}{100}}
\put(2262,667){\blacken\ellipse{100}{100}}
\texture{44555555 55aaaaaa aa555555 55aaaaaa aa555555 55aaaaaa aa555555 55aaaaaa 
	aa555555 55aaaaaa aa555555 55aaaaaa aa555555 55aaaaaa aa555555 55aaaaaa 
	aa555555 55aaaaaa aa555555 55aaaaaa aa555555 55aaaaaa aa555555 55aaaaaa 
	aa555555 55aaaaaa aa555555 55aaaaaa aa555555 55aaaaaa aa555555 55aaaaaa }
\put(1352,1577){\shade\ellipse{150}{150}}
\path(12,2412)(3012,2412)
\path(12,2412)(3012,2412)
\path(312,2712)(312,12)
\path(312,2712)(312,12)
\path(87,1737)(89,1739)(92,1742)
	(99,1748)(109,1758)(123,1772)
	(141,1790)(164,1813)(191,1839)
	(221,1868)(254,1900)(289,1934)
	(326,1970)(363,2006)(401,2042)
	(438,2077)(474,2112)(509,2145)
	(543,2177)(575,2207)(606,2235)
	(635,2261)(663,2286)(689,2309)
	(714,2331)(737,2351)(760,2370)
	(781,2388)(802,2405)(823,2420)
	(843,2435)(862,2449)(880,2462)
	(898,2475)(916,2487)(933,2498)
	(951,2509)(969,2519)(986,2528)
	(1003,2537)(1021,2545)(1038,2553)
	(1055,2559)(1071,2564)(1088,2569)
	(1104,2572)(1120,2575)(1135,2576)
	(1150,2576)(1165,2575)(1179,2573)
	(1192,2570)(1205,2565)(1217,2559)
	(1229,2551)(1240,2543)(1250,2533)
	(1260,2522)(1268,2509)(1277,2496)
	(1284,2481)(1291,2465)(1297,2447)
	(1303,2429)(1308,2409)(1312,2387)
	(1316,2363)(1319,2336)(1322,2308)
	(1325,2279)(1327,2248)(1328,2215)
	(1329,2180)(1330,2144)(1331,2107)
	(1332,2068)(1332,2029)(1333,1989)
	(1333,1947)(1334,1906)(1334,1864)
	(1335,1822)(1336,1781)(1337,1740)
	(1339,1699)(1341,1660)(1343,1622)
	(1346,1585)(1350,1549)(1354,1515)
	(1358,1483)(1364,1452)(1370,1424)
	(1377,1397)(1384,1373)(1392,1351)
	(1402,1330)(1412,1312)(1425,1294)
	(1440,1278)(1456,1266)(1474,1256)
	(1494,1250)(1517,1246)(1542,1246)
	(1570,1248)(1601,1253)(1634,1261)
	(1671,1272)(1710,1285)(1752,1301)
	(1797,1319)(1843,1339)(1889,1361)
	(1936,1383)(1982,1405)(2025,1427)
	(2065,1447)(2099,1465)(2128,1480)
	(2151,1492)(2168,1501)(2178,1507)
	(2184,1511)(2187,1512)
\path(1062,1512)(1065,1513)(1071,1515)
	(1083,1518)(1100,1523)(1124,1530)
	(1153,1538)(1186,1547)(1223,1556)
	(1261,1566)(1300,1575)(1339,1584)
	(1376,1592)(1411,1599)(1444,1604)
	(1474,1608)(1501,1611)(1527,1612)
	(1549,1611)(1570,1610)(1589,1606)
	(1606,1601)(1622,1595)(1637,1587)
	(1651,1578)(1664,1567)(1676,1555)
	(1688,1541)(1698,1526)(1709,1509)
	(1718,1492)(1726,1472)(1734,1452)
	(1741,1430)(1747,1407)(1752,1384)
	(1756,1360)(1760,1335)(1762,1311)
	(1764,1286)(1764,1261)(1764,1237)
	(1763,1213)(1762,1189)(1760,1166)
	(1757,1143)(1754,1121)(1750,1099)
	(1744,1076)(1739,1053)(1733,1030)
	(1726,1008)(1719,985)(1713,962)
	(1706,940)(1699,918)(1693,896)
	(1688,875)(1683,854)(1679,835)
	(1676,816)(1675,798)(1675,782)
	(1677,766)(1680,752)(1685,739)
	(1692,728)(1701,717)(1712,708)
	(1725,699)(1737,693)(1752,688)
	(1769,683)(1788,679)(1810,675)
	(1835,673)(1862,670)(1893,669)
	(1928,668)(1966,668)(2007,668)
	(2052,668)(2099,669)(2148,671)
	(2198,673)(2248,675)(2296,677)
	(2341,679)(2382,681)(2416,683)
	(2443,684)(2463,686)(2476,686)
	(2484,687)(2487,687)
\put(387,237){\makebox(0,0)[lb]{{$\ell_1$}}}
\put(2937,2187){\makebox(0,0)[lb]{{$\ell_2$}}}
\put(387,1737){\makebox(0,0)[lb]{{$r$}}}
\put(912,2262){\makebox(0,0)[lb]{{$q$}}}
\put(1137,1362){\makebox(0,0)[lb]{{$z$}}}
\put(1512,462){\makebox(0,0)[lb]{{$p_1$}}}
\put(2262,462){\makebox(0,0)[lb]{{$p_2$}}}
\put(1962,762){\makebox(0,0)[lb]{{$C_2$}}}
\put(2262,1362){\makebox(0,0)[lb]{{$C_1$}}}
\end{picture}
}

For every choice of splitting $d_1+d_2 = d$ we have ${3d-4}\choose {3d_1-1}$ ways to choose the set of $3d_1-1$ points on $C_1$ from the $3d-4$ points $p_3,\ldots p_{3d-2}$. We have $N_{d_1}$ choices for the curve $C_1$ and $N_{d_2}$ choices for $C_2$. We have $d_1\cdot d_2$ choices for $z$, $d_1$ choices for $q$ and $d_1$ for $r$. This gives

\begin{align*}\deg_B \lambda^{-1}(p_1,&p_2|q,r)\\
&=  \ \  N_d \ \ + 
 \sum_{\begin{array}{c}d=d_1+d_2\\d_1, d_2>0\end{array}} \  {{3d-4}\choose {3d_1-1}}\  \cdot\    N_{d_1}\ N_{d_2}\   \cdot\    d_1 d_2\   \cdot\    d_1^2
\end{align*}

Equating the two sides and rearranging we get the formula.\qed

\subsection{The space of  stable maps}
Gromov--Witten theory allows one to systematically carry out the argument in general, without sweeping things under the rug as I have done above.

Kontsevich introduced the moduli space $\ocM_{g,n}(X,\beta)$ of stable maps, a basic tool in Gromov--Witten theory. As it turned out later, it is a useful moduli space for other purposes.

Fixing a complex projective variety $X$, two integers $g, n\geq 0$ and $\beta \in H_2(X,\ZZ)$, one defines:
$$\ocM_{g,n}(X,\beta) = \left\{(f:C\to X,p_1,\ldots,p_n\in C)\right\},$$
where 
\begin{itemize}
\item $C$ is a nodal connected projective curve,
\item $f:C\to X$ is a morphism with $f_*[C] = \beta$
\item $p_i$ are $n$ distinct points on the smooth locus of $C$, and 
\item The group of automorphisms of $f$ fixing all the $p_i$ is finite.
\end{itemize}

Here an automorphism of $f$ means an automorphism $\sigma:C \to C$ such that $f = f\circ \sigma$, namely a commutative diagram: 

$$\xymatrix{
C\ar[d]_\sigma\ar[rd]^f\\
C\ar[r]_f & X.
}$$

\begin{remark} if $X = $ a point, then $\ocM_{g,n}(X,0) = \ocM_{g,n}$, the Deligne--Mumford stack of stable curves.
\end{remark}

\begin{remark}
It is not too difficult to see that the stability condition on finiteness of automorphisms  is equivalent to either of the following:
\begin{itemize}
\item The sheaf $\omega_C(p_1+\cdots + p_n) \otimes f^*M$ is ample for any sufficiently ample sheaf $M$ on $C$, or
\item We say that a point on the normalization of $C$ is special if it is either a marked point or lies over a node of $C$. The condition is that any rational component $C_0$ of the normalization of $C$ such that $f(C_0)$ is a point, has at least 3 special points, and any such elliptic component has at least one special point.
\end{itemize}
\end{remark}
  
The basic result, treated among other places in \cite{Kontsevich-Manin}, \cite{Fulton-Pandharipande}, is

\begin{theorem}
$\ocM_{g,n}(X,\beta)$ is a proper Deligne--Mumford stack with projective coarse moduli space.
\end{theorem}

\subsection{Natural maps}

The moduli spaces come with a rich structure of maps tying them, and $X$, together. 

\subsubsection{Evaluation}
First, for any $1 \leq i \leq n$ we have natural morphisms, called {\em evaluation morphisms} 

\begin{align*}
\ocM_{g,n}(X,\beta)\ \ \ \ \ &\stackrel{e_i}{\lrar}\ \ \ \ \ X \\ 
(C\stackrel{f}{\to} X,p_1,\ldots,p_n)& \mapsto\ \ \ \ \ \ f(p_i)
\end{align*}

\subsubsection{Contraction}
Next, given a morphism $\phi: X \to Y$ and $n> m$ we get an induced  morphism

\begin{align*}
\ocM_{g,n}(X,\beta) \ \ \ \ \ & \lrar \ \ \ \ \ \ocM_{g,m}(Y,\phi_*\beta)\\
(C\stackrel{f}{\to} X,p_1,\ldots,p_n)& \mapsto\ \mbox{ stabilization of } (C\stackrel{\phi\circ f}{\to} Y,p_1,\ldots,p_m).
\end{align*}
Here in the stabilization we contract those rational components of $C$ which are mapped to a point by $\phi\circ f$ and  have fewer than 3 special points. This is well defined if either $\phi_*\beta\neq 0$ or $2g-2+n>0$.

For instance, if $n>4$ we get a morphism
$\ocM_{0,n}(X,\beta) \to \ocM_{0,4}$.

\subsection{Boundary of moduli} Understanding the subspace of maps with degenerate source curve $C$ is key to Gromov--Witten theory.

\subsubsection{Fixed degenerate curve}
Suppose we have a degenerate curve \begin{center}
\setlength{\unitlength}{0.00083333in}
\begingroup\makeatletter\ifx\SetFigFont\undefined%
\gdef\SetFigFont#1#2#3#4#5{%
  \reset@font\fontsize{#1}{#2pt}%
  \fontfamily{#3}\fontseries{#4}\fontshape{#5}%
  \selectfont}%
\fi\endgroup%
{\renewcommand{\dashlinestretch}{30}
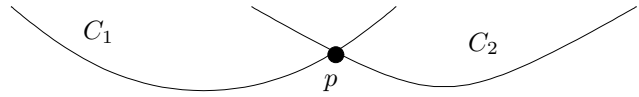
\begin{picture}(3924,610)(0,-10)
\put(2037,283){\blacken\ellipse{100}{100}}
\path(12,583)(13,582)(15,581)
	(18,578)(23,574)(30,567)
	(40,559)(53,549)(68,536)
	(86,522)(106,505)(129,487)
	(155,467)(183,446)(212,424)
	(244,401)(277,377)(312,353)
	(348,329)(385,306)(424,282)
	(464,259)(505,237)(547,215)
	(591,194)(637,174)(684,155)
	(734,137)(786,121)(839,106)
	(896,92)(955,81)(1016,71)
	(1080,64)(1146,60)(1212,58)
	(1278,60)(1344,64)(1408,71)
	(1469,81)(1528,92)(1585,106)
	(1638,121)(1690,137)(1740,155)
	(1787,174)(1833,194)(1877,215)
	(1919,237)(1960,259)(2000,282)
	(2039,306)(2076,329)(2112,353)
	(2147,377)(2180,401)(2212,424)
	(2241,446)(2269,467)(2295,487)
	(2318,505)(2338,522)(2356,536)
	(2371,549)(2384,559)(2394,567)
	(2401,574)(2406,578)(2409,581)
	(2411,582)(2412,583)
\path(1512,583)(1514,582)(1517,580)
	(1524,576)(1535,570)(1550,561)
	(1571,549)(1596,535)(1626,518)
	(1661,498)(1700,477)(1743,453)
	(1788,428)(1836,403)(1885,376)
	(1935,350)(1985,324)(2035,299)
	(2084,275)(2132,252)(2178,230)
	(2223,210)(2266,191)(2307,174)
	(2347,159)(2385,145)(2422,132)
	(2457,121)(2491,112)(2525,104)
	(2557,97)(2589,92)(2620,88)
	(2651,85)(2681,84)(2712,83)
	(2743,84)(2773,85)(2804,88)
	(2835,92)(2867,97)(2899,104)
	(2933,112)(2967,121)(3002,132)
	(3039,145)(3077,159)(3117,174)
	(3158,191)(3201,210)(3246,230)
	(3292,252)(3340,275)(3389,299)
	(3439,324)(3489,350)(3539,376)
	(3588,403)(3636,428)(3681,453)
	(3724,477)(3763,498)(3798,518)
	(3828,535)(3853,549)(3874,561)
	(3889,570)(3900,576)(3907,580)
	(3910,582)(3912,583)
\put(462,358){\makebox(0,0)[lb]{{$C_1$}}}
\put(2862,283){\makebox(0,0)[lb]{{$C_2$}}}
\put(1962,58){\makebox(0,0)[lb]{{$p$}}}
\end{picture}
}

\end{center}
$$C = C_1 \mathop{\cup}\limits^p C_2.$$
So $C$ is a fibered coproduct of two curves.
By the universal property of coproducts
\begin{align*} 
Hom(C, X) 
&= Hom(C_1,X) \mathop\times\limits_{Hom(p, X)}Hom(C_2,X) \\  
&= Hom(C_1,X)\ \ \  \mathop\times\limits_{X}\ \ \ Hom(C_2,X) 
\end{align*}

\subsubsection{Varying degenerate curve: the boundary of moduli}
We can work this out in the fibers of universal the families. If we set $g = g_1+g_2,$ $n= n_1+n_2$ and $\beta = \beta_1 + \beta_2$ we get a morphism
$$ \ocM_{g_1, n_1+1}(X, \beta_1) \times_X  \ocM_{g_2, n_2+1}(X, \beta_2) \lrar \ocM_{g, n}(X, \beta), $$ with the fibered product over $e_{n_1+1}$ on the left and $e_{n_2+1}$ on the right. On the level of points this is obtained by gluing curves $C_1$ at point $n_1+1$ with $C_2$ at point $n_2+1$ and matching the maps $f_1, f_2$. This is a finite unramified map, and we can think of the product on the left as a space of stable  $n$-pointed maps of genus $g$ and class $\beta$ with a distinguished marked node.

\subsubsection{Compatibility with evaluation maps}
This map is automatically compatible with the other ``unused" evaluation maps.
For instance if $i\leq n_1$ we get a commutative diagram
$$\xymatrix{ \ocM_{g_1, n_1+1}(X, \beta_1) \mathop{\times}\limits_X  \ocM_{g_2, n_2+1}(X, \beta_2) \ar[r]\ar_{\pi_1}[d] &\ocM_{g, n}(X, \beta)\ar^{e_i}[d]\\
\ocM_{g_1, n_1+1}(X, \beta_1)\ar^{e_i}[r]& X.
} $$

As these tend to get complicated we will give the marking their ``individual labels" rather than a number (which may change through the gluing or contraction maps).

\subsection{Gromov--Witten classes}
We start with some simplifying assumptions:
\begin{itemize}
\item $g=0$
\item $X$ is ``convex" \cite{Fulton-Pandharipande}: no map of rational curve to $X$ is obstructed. Examples: $X = \PP^r$, or any homogeneous space.
\end{itemize}

We simplify notation: $M = \ocM_{0,n+1}(X, \beta)$. We take $\gamma_i \in H^*(X, \QQ)^{even}$ to avoid sign issues.

We now define {\em Gromov--Witten classes}:
\begin{definition}
$$\<\gamma_1,\ldots,\gamma_n,*\>^X_\beta:= e_{n+1\,*}(e^*(\gamma_1\times\ldots\times\gamma_n)\cap[M]))\in H^*(X).$$ Here the notation is $e:=e_1\times \cdots\times e_n : M\to X^n$: 
$$\xymatrix{ M\ar^{e_{n+1}}[rr]\ar^{e:=e_1\times \cdots\times e_n}[dd] && X \\ \\ X^n. 
}$$
\end{definition}

When $X$ is fixed we will suppress the superscript $X$ from the notation.

\subsection{The WDVV equations}
The main formula in genus 0 Gromov--Witten theory is the Witten--Dijkgraaf--Verlinde--Verlinde (or WDVV) formula. It is simplest to state for $n=3$:

\begin{theorem}
$$
\sum_{\beta_1+\beta_2=\beta}\Big\<\<\gamma_1,\gamma_2,*\>_{\beta_1},\gamma_3,* \Big\>_{\beta_2}=
\sum_{\beta_1+\beta_2=\beta}\Big\<\<\gamma_1,\gamma_3,*\>_{\beta_1},\gamma_2,* \Big\>_{\beta_2}.
$$

For general $n\geq 3$ it is convenient to label the markings by a finite set $I$ to avoid confusion with numbering. The formula is

\begin{align*}
&\sum_{\beta_1+\beta_2=\beta} \sum_{A\sqcup B = I}
\Big\<\<\gamma_1,\gamma_2,\, \delta_{A_1},\ldots,\delta_{A_k},\, *\>_{\beta_1},\gamma_3,\, 
\delta_{B_1},\ldots,\delta_{B_m},\, * \Big\>_{\beta_2} \\
=&\sum_{\beta_1+\beta_2=\beta}\sum_{A\sqcup B = I}
\Big\<\<\gamma_1,\gamma_3,\, \delta_{A_1},\ldots,\delta_{A_k},\, *\>_{\beta_1},\gamma_2,\, 
\delta_{B_1},\ldots,\delta_{B_m},\, * \Big\>_{\beta_2}
\end{align*}
\end{theorem}
Note that the only thing changed between the last two lines is the placement of $\gamma_2$ and $\gamma_3$.

This formalism of Gromov--Witten classes and the WDVV equations is taken from a yet unpublished paper of Graber and Pandharipande. One advantage is that it works with cohomology replaced by Chow groups. 

\subsubsection{Gromov--Witten numbers}
The formalism that came to us from the physics world is equivalent, though different, and involves Gromov--Witten numbers, defined as follows:

\begin{definition}
$$\<\gamma_1,\ldots,\gamma_n\>_\beta^X\ :=\ \int_Me^*(\gamma_1\times\ldots\times\gamma_n),$$
where this time $M = \ocM_{0,n}(X, \beta)$. 
\end{definition}

The following elementary lemma allows one to go back and forth between these formalisms:
\begin{lemma}\hfill
\begin{enumerate}
\item $\displaystyle \<\gamma_1,\ldots,\gamma_n,\gamma_{n+1}\>_\beta\ =\ \int_X \<\gamma_1,\ldots,\gamma_n,*\>_\beta\, \cup\, \gamma_{n+1}$
\item Choose a basis $\{\alpha_i\}$ for $H^*(X , \QQ)$. Write the intersection matrix $\int_X \alpha_i \cup \alpha_j = g_{ij}$, and denote the inverse matrix entries by $g^{ij}$. Then
$$\<\gamma_1,\ldots,\gamma_n,*\>_\beta\ \   =\ \  \sum_{i,j}\ \<\gamma_1,\ldots,\gamma_n,\alpha_i\>_\beta\, g^{ij}\ \alpha_j.$$ 
\end{enumerate}
\end{lemma}

The WDVV equations then take the form
\begin{align*}&\sum_{\beta_1+\beta_2=\beta}\ \ \sum_{A\sqcup B=I}\ \ \sum_{i,j} \ \ 
\<\gamma_1,\gamma_2,\delta_A,\alpha_i\>_{\beta_1}\ \ g^{ij}\ \ \<\alpha_j,\gamma_3,\delta_B,\gamma_4\>_{\beta_2}\\
=&\sum_{\beta_1+\beta_2=\beta}\ \ \sum_{A\sqcup B=I}\ \ \sum_{i,j} \ \ 
\<\gamma_1,\gamma_3,\delta_A,\alpha_i\>_{\beta_1}\ \ g^{ij}\ \ \<\alpha_j,\gamma_2,\delta_B,\gamma_4\>_{\beta_2}
\end{align*}

\subsection{Proof of WDVV}
We sketch the proof of WDVV, in the formalism of  Gromov--Witten classes,  under strong assumptions:
\begin{itemize}
\item $e_i: M \lrar X$ is smooth, and
\item the contraction $M \lrar \ocM_{0,4}$ is smooth, and moreover each node can be smoothed  out independently. 
\end{itemize}

These assumptions hold for the so called convex varieties discussed by Fulton and Pandharipande.

\subsubsection{Setup} Now, it would be terribly confusing to use the usual numbering for the markings, the evaluation maps  and the cohomology classes pulled back by the corresponding evaluation marks, as we will use the structure of the boundary discussed above. Instead we give them names.

The first moduli space we need is the space of genus 0 pointed stable maps to X with class $\beta_1$. The markings are used in three ways:
\begin{enumerate}
\item The first two are used to pull back $\gamma_1,\gamma_2$ on the left hand side (and $\gamma_1,\gamma_3$ on the right).
\item The next bunch is used to pull back the $\delta_{A_i}$.
\item The last is used to push forward.
\end{enumerate}
It is convenient to put together the first two sets of and call the result $\hat A$. The last marking will be denoted by the symbol $\blacktriangleright$, suggesting something is to be glued on the right. For short notation we will use
$$\ocM_1 = \ocM_{0, \hat A \sqcup \blacktriangleright}(X,\beta_1).$$
We  also use shorthand for the cohomology and homology classes:

\begin{align*}
\eta_1\ \  &=\ \  e_{\hat A}^*(\gamma_1 \times \gamma_2 \times \delta_A)\\
\xi_1\ \  &=\  \ \eta_1\  \cap\  \left[ M_1 \right]
\end{align*}

The second moduli space we need is the space of genus 0 pointed stable maps to X with class $\beta_2$. The markings are used in four ways:
\begin{enumerate}
\item The first will be used to pull back $\<\gamma_1,\gamma_2,\, \delta_{A_1},\ldots,\delta_{A_k},\, *\>_{\beta_1}$. It is the marking used for gluing and is accordingly denoted $\blacktriangleleft$.
\item The next one is used to pull back $\gamma_3$ on the left hand side (and $\gamma_2$ on the right).
\item The next bunch is used to pull back the $\delta_{B_i}$.
\item The last is used to push forward.
\end{enumerate}
We leave the first one alone, we put together the next two sets of and call the result $\hat B$. The last marking will be denoted by the symbol $\bullet$. The notation we will use is
$$\ocM_2 = \ocM_{0,\blacktriangleleft\sqcup \hat B \sqcup \bullet}(X,\beta_1).$$
The shorthand for the cohomology and homology classes is:

\begin{align*}
\eta_2\ \  &=\ \  e_{\hat B}^*(\gamma_3 \times \delta_B)\\
\xi_1\ \  &=\  \ \eta_2\  \cap\  \left[ M_2 \right]
\end{align*}

The key to the proof of WDVV is to express the class on either side in terms of  something symmetric on the glued moduli space. We use the notation
$$ \ocM_1\times_X \ocM_1 \lrar \ocM := \ocM_{0, \hat A \sqcup \hat B \sqcup \bullet}(X, \beta),$$ where the gluing is done with respect to $e_\blacktriangleright: \ocM_1 \to X$ on the left, and $e_\blacktriangleleft: \ocM_2 \to X$ on the right. 
A relevant symmetric class  which comes up as a bridge between the two sides of the formula is
$$\eta_{12}\ \  = \ \ e_{\hat A \sqcup \hat B}^*\big( \gamma_1 \times \gamma_2 \times \delta_A \times \gamma_3 \times \delta_B\big).$$

\subsubsection{The fibered product diagram.}

Consider the diagram

$$\xymatrix{
\ocM_1 \times_X\ocM_2\ar^{p_2}[r]\ar_{p_1}[d] & \ocM_2\ar^{e_\bullet}[r] \ar^{e_\blacktriangleleft}[d] & X \\
\ocM_1 \ar^{e_\blacktriangleright}[r] & X.
}$$
It is important that $e_\blacktriangleleft$ is smooth, in particular flat. We suppress Poincar\'e duality isomorphisms in the notation. 

We now have by definition
\begin{align*}\Big\<\<\gamma_1,\gamma_2,\, \delta_{A},\, *\>_{\beta_1},\gamma_3,\, 
\delta_{B},\, * \Big\>_{\beta_2}& = 
e_{\bullet\,*}\Big(e_\blacktriangleleft^*\left(e_{\blacktriangleright\,*}\xi_1\right)\cap \xi_2\Big) \\
\intertext{and the projection formula gives}
& = (e_\bullet\circ p_2)_*\Big(p_1^*\eta_1 \cup p_2^* \eta_2\cap [\ocM_1\times_X\ocM_2]\Big)
\end{align*}

So we need to understand the class $p_1^*\eta_1 \cup p_2^* \eta_2\cap [\ocM_1\times_X\ocM_2]$ on the fibered product.

\subsubsection{End of proof.}

Since we have a big sum in the WDVV equation, we need to take the union of  all these fibered products. They are put together by the gluing maps, as in the following fiber diagram:
$$\xymatrix{
\coprod\limits_{\beta_1+\beta_2 = \beta} \coprod\limits_{A\sqcup B = I} \ocM_1\times_X\ocM_2 \ar[r]\ar@/^1.5pc/[rr]^{\ell}
 & \Delta^X_{(12|3\bullet)} \ar@{^(->}[r]\ar[d] & \ocM \ar^{st}[d]\\
 &\{(12|3\bullet)\}\ar@{^(->}[r] & \ocM_{0,4}
}$$
where the markings in $\ocM_{0,4}$ are denoted $1,2,3,$ and $\bullet$ to match with our other notation. The subscheme $\Delta^X_{(12|3\bullet)}\subset \ocM$ is defined by this diagram.

The smoothness assumption on the stabilization/contraction map $$st: \ocM \to \ocM_{0,4}$$ guarantees that the gluing map $$j: \coprod\limits_{\beta_1+\beta_2 = \beta} \coprod\limits_{A\sqcup B = I} \ocM_1\times_X\ocM_2 \to \Delta^X_{(12|3\bullet)}$$ is finite and birational, and in fact we have an equality in Chow classes: 
$$\ell_* \left[\coprod\limits_{\beta_1+\beta_2 = \beta} \coprod\limits_{A\sqcup B = I} \ocM_1\times_X\ocM_2\right] \ \ = \ \ st^* \left[ \{(12|3\bullet)\}\right].$$

We can now complete our equation:
\begin{align*}\sum\sum\Big\<\<\gamma_1,\gamma_2,\, \delta_{A},\, *\>_{\beta_1}&,\gamma_3,\, 
\delta_{B},\, * \Big\>_{\beta_2}\\
& = 
\sum\sum  (e_\bullet\circ p_2)_*\Big(p_1^*\eta_1 \cup p_2^* \eta_2\cap [\ocM_1\times_X\ocM_2]\Big)\\
& = \sum\sum  (e_\bullet\circ p_2)_*\Big(\ell^* \eta_{12}\cap [\ocM_1\times_X\ocM_2]\Big)\\
& = e_{\bullet \, *}\Big(\eta_{12}\cap st^*[ \{(12|3\bullet)\}]\Big),
\end{align*}
where the last evaluation map is $e_\bullet: \ocM \to \ocM_{0,4}$. 

The latter expression is evidently symmetric and therefore the order of $\gamma_2$ and $\gamma_3$ is immaterial. The formula follows. \qed

\subsection{About the general case}

The WDVV equation holds in general when one replaces fundamental classes of moduli space in the smooth case with virtual fundamental classes, something that takes care of the non-smoothness in an organized manner. The drawback of virtual fundamental classes is that in general one loses the enumerative nature of Gromov--Witten classes. However the flexibility of the formalism allows one to compute cases where the classes are enumerative by going through cases where they are not. 

I will definitely not try to get into the details here - they belong in a different lecture series. However I cannot ignore the subject completely - in the orbifold situation virtual fundamental classes are always necessary!  Let me just indicate the principles.

First, the smoothness assumption of stabilization holds in the ``universal case'' - by which I mean to say that it holds for the moduli stacks of pre-stable curves introduced by Behrend \cite{Behrend}. We in fact have a fiber diagram extending the above:
$$\xymatrix{
\coprod\limits_{\beta_1+\beta_2 = \beta} \coprod\limits_{A\sqcup B = I} \ocM_1\times_X\ocM_2 \ar[r]\ar@/^1.5pc/[rr]^{\ell}\ar[d]
 & \Delta^X_{(12|3\bullet)} \ar@{^(->}[r]\ar[d] & \ocM \ar^{\mbox{forgetful}}[d]\\
 \fM_{03}\times \fM_{0,3} \ar[r] & \fD_{0,4}\ar[r]\ar[d]& \fM_{0,4}\ar^{\mathfrak{st}}[d] \\
 &\{(12|3\bullet)\}\ar@{^(->}[r] & \ocM_{0,4}
}$$
where $\mathfrak{st}: \fM_{0,4} \lrar  \ocM_{0,4}$ is flat. This is a situation where refined pull-backs can be used.

The formalism of algebraic virtual fundamental classes works for an arbitrary smooth projective $X$. It requires the definition of classes $[\ocM]^{\vir}$ in the Chow group of each  moduli space $\ocM$. The key property that this satisfies is summarized as follows:

Consider the diagram
$$\xymatrix{
\ocM_1 \times \ocM_2 \ar[d] & \ocM_1 \times_X \ocM_2 \ar[l]\ar[d]\ar@{^(->}[r] & \sqcup \ocM_1\times_X\ocM_2\ar[r]\ar[d] & \ocM\ar[d]\\
X \times X & X \ar^{\Delta}[l] &  \fM_{03}\times \fM_{0,3} \ar^{\mathfrak{gl}}[r] & \fM_{0,4}.}$$

The condition is:

$$\sum_{\beta_1+\beta_2=\beta}\sum_{A\sqcup B = I}
\Delta^!\left([\ocM_1]^\vir\times [\ocM_2]^\vir\right) \ \ = \ \ \mathfrak{gl}^![\ocM]^\vir. $$


For a complete algebraic treatment and an explanation why this is the necessary  equation see \cite{Behrend-Fantechi}.

The main theorem is 
\begin{theorem}
This equation holds for the class associated to obstruction theory
$$[\ocM]^\vir \ \ = \ \ (\ocM, E)$$
with $$E = \bR\pi_*f^*T_X$$ coming from the diagram of the universal curve
 $$\xymatrix{ \cC \ar_\pi[d] \ar[r]^f & X \\ \ocM
 }$$
\end{theorem}

\section{Orbifolds / Stacks}
\subsection{Geometric orbifolds}
One can spend an entire lecture series laying down the foundations of orbifolds or stacks. Since I do not have the luxury, I'll stick to a somewhat intuitive, and necessarily imprecise, presentation. The big drawback is that people who do not know the subject get only a taste of what it is about, and have to look things up to really understand what is going on. Apart from the standard references, one may consult the appendix of \cite{Vistoli} or \cite{Edidin} for an introduction to the subject.

I will use the two words ``orbifold" and ``stack" almost interchangeably.

Geometrically, an orbifold $\cX$ is locally given as a quotient of a space, or manifold, $Y$ by the action of a finite group, giving a chart $[Y/G] \to \cX$. The key is to remember something about the action.

A good general way to do this is to think of $\cX$ as the equivalence class of a groupoid
$$R \mathop{\double}\limits_s^t V$$
where
\begin{itemize}
\item  $s,t$ are \'etale morphisms or, more generally, smooth morphisms, and
\item  the morphism $s\times t:R \to V$ is required to be finite.
\end{itemize}

The notation used is $\cX \simeq [R \double V]$. You think of $R$ as an equivalence relation on the points of $V$, and of $\cX$ as the ``space of equivalence classes".

This notation of a groupoid and its associated orbifold is very much a shorthand - the complete data requires a composition map 
$$R \mathop\times\limits_{^tV^s}R\ \  \to\ \  R$$
as well as an identity morphism $V \to R$, satisfying standard axioms I will not write in detail. They are inspired by the case of a quotient:

A quotient is the orbifold associate to the following groupoid:  $[V/G] = [R\double V]$, where $R = G \times V$, the source map $s:G\times V$ is the projection, and the target map $t: R \to V$ is the action.

I have not described the notion of equivalence of groupoids, neither did I describe the notion of a morphism of orbifolfds presented by groupoids. It is a rather complicated issue which I had rather avoid. The moduli discussion below will shed some light on it.

\subsection{Moduli stacks}

Orbifolds come about rather frequently in the theory of moduli. In fact, the correct definition of an algebraic stack is as a sort of tautological solution to a moduli problem.

An algebraic stack $\cX$ is by definition a category, and implicitly the category of families of those object we want to parametrize. It comes with a ``structure functor" $\cX \to Sch$ to the category of scheme, which to each family associates the base of the family. A morphism of stacks is a functor commuting with the structure functor.

Here is the key example: $\cM_g$ - the moduli stack of curves. The category $\cM_g$ has as its objects $$\left\{\begin{array}{c}\cC\\ \dar\\B\end{array}\right\}$$ where each $\cC \to B$ is a family of curves of genus $g$. The structure functor $\cM_g \to Sch$ of course sends the object $\cC \to B$ to the base scheme $B$. It is instructive to understand what arrows we need. Of course we want to classify families up to isomorphisms, so isomoprhisms 
$$ \xymatrix{ \cC \ar[d]\ar^{\sim}[r] &\cC'\ar[d] \\ B \ar@{=}[r]& B
}$$
must be included. But it is not hard to fathom that pullbacks are important as well, and indeed a morphism in $\cM_g$ is defined to be a {\em fiber diagram} 
$$ \xymatrix{ \cC \ar[d]\ar[r] &\cC'\ar[d] \\ B \ar[r]& B'
}$$
That's why you'll hear the term ``fibered category" used.

How is a scheme thought of as a special case of a stack? A scheme $X$ is the moduli space of its own points. So an object is ``a family of points of $X$ parametrized by $B$", i.e. a morphism $B \to X$. So the stack associated to $X$ is just the category of schemes over $X$.

Here is another example: consider a quotient orbifold $[Y/G]$, where $Y$ is a variety and $G$ a finite group. How do we think of it as a category? A point of the orbit space $Y/G$ should be an orbit of $G$ in $Y$, and you can think of an orbit as the image of an equivariant map of a principal homogeneous  $G$-space $P$ to $Y$. We make this the definition: 
\begin{itemize}
\item an object of $[Y/G]$ over a base scheme $B$ is a principal homogeneous $G$-space $P \to B$ together with a $G$-equivariant map $P \to Y$:
$$\xymatrix{ P \ar[d]\ar[r]& Y \\B.}$$
\item An arrow is a fiber diagram:
$$\xymatrix{ P \ar[d]\ar[r]\ar@/^.8pc/[rr]&P_1\ar[r]\ar[d] &Y \\B \ar[r] &B_1 .}$$
\end{itemize}

An important special case is when $Y$ is a point, giving the classifying stack of $G$, denoted $$\cB G := [\{pt\}/G].$$ Objects are just principal $G$-bundles, and morphisms are fiber diagrams.

Kai Behrend gave an elegant description of the stack associated to a groupoid in general. First we note the following: for any scheme $X$ and any \'etale surjective $V \to X$, one can write $R_V = V \times_X V$ and the two projections give a groupoid $R_V \double V$. If, as suggested above, $R_V$ is to be considered as an equivalence relation on $V$, then clearly the equivalence classes are just points of $X$, so we had better define things so that $X = [R_V\double V]$.  

Now given a general groupoid $R \double V$, an object over a base scheme $B$ is very much like a principal homogeneous space: it consists of an \'etale covering $U \to B$, giving rise to $R_U \double U$ as above, together with maps $U \to V$ and $R_U \to R$ making the following diagram  (and all its implicit siblings) {\em cartesian}:
$$\xymatrix{
R_U\ar[r]\doubledown & R\doubledown\\
U\ar[r]\ar[d] & V\\
B
}$$

There is an important object of $\cX = [R \double V]$ with the scheme $V$ as its base: you take $U = R$ above, with the two maps $U\to B$ and $U\to V$ being the source and target maps $R\to V$ respectively. What it does is it gives an {\em \'etale covering $V \to \cX$}. The existence of such a thing is in fact an axiom required of a fibered category  to be a Deligne--Mumford algebraic stack, but since I have not gotten into details you'll need to study this elsewhere. The requirement says in essence that every object should have a universal deformation space.

Some words you will see: 
\begin{itemize}
\item an algebraic space is a stack of the form $[R\double V]$ with $R \to V \times V$ injective. This is where ``stacks" meet ``sheaves".
\item An Artin stack, also known as a  general algebraic stack, is what you get when you only require the source and target maps $R\to V$ to be smooth, and do not require $ R \to V \times V$ to be proper either. You can't quite think of an Artin stack as ``locally the quotient of a scheme by a group action" - the categorical viewpoint is necessary.
\end{itemize}

\subsection{Where do stacks come up?}
\subsubsection{Moduli, of course}
The first place where you meet stacks is when trying to build moduli spaces. Commonly, fine moduli spaces do not exist because objects have automorphisms, and stacks are the right replacement.
\subsubsection{Hidden smoothness}
But eve if you are not too excited by moduli spaces, stacks, in their incarnation as orbifolds,  are here to stay. The reason is, it is often desirable to view varieties with finite quotient singularities as if it were smooth, and indeed to every such variety there is  a relevant stack, which is indeed smooth. 

This feature comes up, and increasing in appearance, in many topics in geometry: the minimal model program, mirror symmetry, geometry of 3-manifolds, the McKay correspondence, and even in Haiman's $n!$ theorem (though this is not the way Haiman would present it).

\subsection{Attributes of orbifolds}
If one is to study orbifolds along with varieties or manifolds, one would like to have tools similar to ones available for varieties and manifolds.

Indeed, the theory is well developed.
\begin{itemize}
\item Orbifolds have homology and cohomology groups, and cohomology of smooth orbifolds satisfies Poincar\'e duality  with rational coefficients. 
\item Chow groups with rational coefficients for Deligne--Mumford stacks  were constructed by Vistoli and Gillet independently in the 80s. More recently Kresch showed in his thesis that they have Chow groups with integer coefficients. 
\item One can talk about sheaves, K-theory and derived categories of stacks. That's a natural framework for the McKay correspondence.
\item Smooth Deligne--Mumford stack have a dualizing invertible sheaf.
\item Laumon and Moret-Bailly introduced the cotangent complex $\LL_{\cX}$ of an algebraic stack. This is rather easy for Deligne--Mumford stacks, but their construction was found to be flawed for Artin stack. This problem was recently corrected by Martin Olsson \cite{Olsson-stacks}.
\item Deligne--Mumford stacks have coarse moduli spaces - this is a theorem of Keel and Mori \cite{Keel-Mori}.  For $[Y/G]$ the moduli space is just the geometric quotient $Y/G$, namely the orbit space. In general this is an algebraic space $X$ with a morphism $\cX \to X$ which is universal, and moreover such that $\cX(k)/\Isom \to X(k)$ is bijective whenever $k$ is an algebraically closed  field.
\item The inertia stack: this is a natural stack associated to $\cX$, which in a way points to where $\cX$ fails to be a space. Every object $\xi$ of $\cX$ has its automorphism group $\Aut(\xi)$, and these can be put together in one stack $\cI(\cX)$, whose objects are pairs $(\xi,\sigma)$, with $\xi$ an object of $\cX$ and $\sigma\in \Aut(\xi)$. One needs to know that this is an algebraic stack when $\cX$ is, a Deligne--Mumford stack when $\cX$ is, etc. This follows from the abstract and not - too - illuminating formula
$$\cI(\cX) = \cX \mathop\times\limits_{\cX \times \cX} \cX,$$
where the product is taken relative to the diagonal map on both sides. As an example, if $\cX = \cB G$ then $\cI(\cX) = [G/G]$, where $G$ acts on itself by conjugation.
\end{itemize}

There is one feature which one likes to ignore, but there comes a point where one needs to face the facts of life: stacks are not a category! They are a 2-category: arrows are functors, 2-arrows are natural transformations.

\subsection{\'Etale gerbes}

This is a class of stacks which will come up in our constructions. 

Informally, an \'etale gerbe is a stack which locally (in the \'etale topology) looks like $X \times \cB G$, with $G$ a finite group.

Formally (but maybe not so intuitively), it is a Deligne--Mumford  stack $\cX $ such that the morphisms $\cI(\cX) \to \cX \to X$ are all finite \'etale.

We need to tie this thing better with a group $G$. 
We will only need to consider the case where $G$ is abelian.

An \'etale gerbe $\cG \to X$ is said to be {\em banded}  by the finite abelian group $G$ if one is given, for every object $\xi\in \cG(S)$, an isomorphism $G(S) \simeq \Aut(\xi)$ in a functorial manner.

Gerbes banded by $G$ can be thought of as principal homogeneous spaces under the ``group stack" $\cB G$. As such, they are classified by the ``next cohomology group over", $H^2_{\text{\'et}}(X,G)$. This led Giraud in his thesis under Grothendieck to define the non-abelian second cohomology groups using non-abelian gerbes, but this goes too far afield for us.

\section{Twisted stable maps}
\subsection{Stable maps to a stack}

Consider a semistable elliptic surface, with base $B$ and a section. We can naturally view this as a map $B \to \ocM_{1,1}$. Angelo Vistoli, when he was on sabbatical at Harvard in 1996, asked the following beautiful, and to me very inspiring, question: what's a good way to compactify the moduli of elliptic surfaces? can one use stable maps to get a good compactification?

Now consider in general:
$$\xymatrix{\cX \ar[d] & \text{ Deligne--Mumford stack with}\\ X & \text{ projective coarse moduli space}
}$$
In analogy to $\ocM_{g,n}(X, \beta)$, we want a compact moduli space of maps $C \to \cX$.

One can define stable maps as in the scheme case, but there is a problem:  the result is not compact. As Angelo Vistoli likes to put it, trying to work with a non-compact moduli space is like trying to keep your coins when you have holes in your pockets. The solution that comes naturally is that 
\begin{quote} \bf the source curve $\cC$ must acquire a stack structure as well as it degenerates! 
\end{quote}

Both problem and solution are clearly present in the following example, which is ``universal" in the sense that we take $\cX$ to be a one parameter family of curves itself:

Consider $\PP^1\times \PP^1$ with coordinates $x,s$ near the origin and the projection with coordinate $s$  onto $\PP^1$. Blowing up the origin we get a family of curves, with general fiber $\PP^1$ and special fiber a nodal curve, with local equation $xy = t$ at the node. Taking base change $\PP^1 \to \PP^1$ of degree 2 with equation $t^2 = s$ we get a singular scheme $X$ with a map $X \to \PP^1$ given by coordinate $s$. This is again a family of $\PP^1$'s with nodal special fiber, but local equation $xy = s^2$.

This is a quotient singularity, and using the chart $[\AA^2/ (\ZZ/2\ZZ)]$ with coordinates $u,v$ satisfying  $u^2 = x, v^2 = y$ we get a smooth orbifold $\cX$, with coarse moduli space $X$ and   a map $\cX \to \PP^1$. It is a family of $\PP^1$'s parametrized by $\PP^1$, degenerating to an orbifold curve.

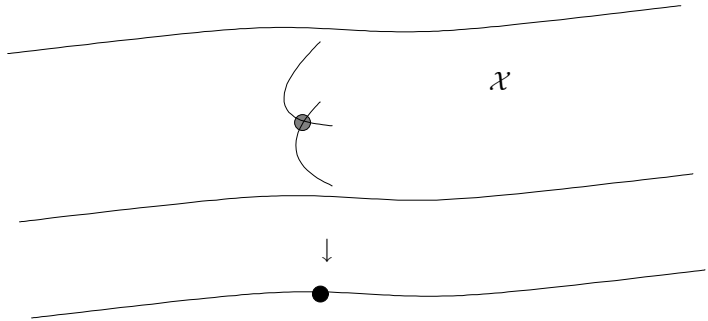
\setlength{\unitlength}{0.00083333in}
\begingroup\makeatletter\ifx\SetFigFont\undefined%
\gdef\SetFigFont#1#2#3#4#5{%
  \reset@font\fontsize{#1}{#2pt}%
  \fontfamily{#3}\fontseries{#4}\fontshape{#5}%
  \selectfont}%
\fi\endgroup%
{\renewcommand{\dashlinestretch}{30}
\begin{picture}(4374,1989)(0,-10)
\put(1962,162){\blacken\ellipse{100}{100}}
\texture{44555555 55aaaaaa aa555555 55aaaaaa aa555555 55aaaaaa aa555555 55aaaaaa 
	aa555555 55aaaaaa aa555555 55aaaaaa aa555555 55aaaaaa aa555555 55aaaaaa 
	aa555555 55aaaaaa aa555555 55aaaaaa aa555555 55aaaaaa aa555555 55aaaaaa 
	aa555555 55aaaaaa aa555555 55aaaaaa aa555555 55aaaaaa aa555555 55aaaaaa }
\put(1850,1232){\shade\ellipse{100}{100}}
\path(12,1662)(14,1662)(18,1663)
	(25,1664)(37,1665)(54,1667)
	(76,1670)(104,1673)(138,1677)
	(178,1682)(223,1687)(273,1693)
	(328,1699)(385,1706)(446,1713)
	(508,1720)(572,1728)(636,1735)
	(699,1742)(762,1749)(824,1756)
	(883,1763)(941,1769)(997,1775)
	(1050,1780)(1101,1786)(1150,1791)
	(1196,1795)(1240,1799)(1282,1803)
	(1322,1807)(1360,1810)(1396,1813)
	(1431,1815)(1464,1817)(1495,1819)
	(1526,1821)(1556,1822)(1584,1824)
	(1612,1824)(1648,1825)(1683,1826)
	(1718,1826)(1752,1826)(1785,1826)
	(1818,1825)(1851,1825)(1884,1824)
	(1916,1822)(1949,1821)(1982,1819)
	(2014,1818)(2047,1816)(2079,1814)
	(2112,1812)(2145,1810)(2177,1808)
	(2210,1806)(2242,1805)(2275,1803)
	(2308,1802)(2340,1800)(2373,1799)
	(2406,1799)(2439,1798)(2472,1798)
	(2506,1798)(2541,1798)(2576,1799)
	(2612,1799)(2640,1800)(2668,1802)
	(2698,1803)(2729,1805)(2760,1807)
	(2793,1809)(2828,1811)(2864,1814)
	(2902,1817)(2942,1821)(2984,1825)
	(3028,1829)(3074,1833)(3123,1838)
	(3174,1844)(3227,1849)(3283,1855)
	(3341,1861)(3400,1868)(3462,1875)
	(3525,1882)(3588,1889)(3652,1896)
	(3716,1904)(3778,1911)(3839,1918)
	(3896,1925)(3951,1931)(4001,1937)
	(4046,1942)(4086,1947)(4120,1951)
	(4148,1954)(4170,1957)(4187,1959)
	(4199,1960)(4206,1961)(4210,1962)(4212,1962)
\path(87,612)(89,612)(93,613)
	(100,614)(112,615)(129,617)
	(151,620)(179,623)(213,627)
	(253,632)(298,637)(348,643)
	(403,649)(460,656)(521,663)
	(583,670)(647,678)(711,685)
	(774,692)(837,699)(899,706)
	(958,713)(1016,719)(1072,725)
	(1125,730)(1176,736)(1225,741)
	(1271,745)(1315,749)(1357,753)
	(1397,757)(1435,760)(1471,763)
	(1506,765)(1539,767)(1570,769)
	(1601,771)(1631,772)(1659,774)
	(1687,774)(1723,775)(1758,776)
	(1793,776)(1827,776)(1860,776)
	(1893,775)(1926,775)(1959,774)
	(1991,772)(2024,771)(2057,769)
	(2089,768)(2122,766)(2154,764)
	(2187,762)(2220,760)(2252,758)
	(2285,756)(2317,755)(2350,753)
	(2383,752)(2415,750)(2448,749)
	(2481,749)(2514,748)(2547,748)
	(2581,748)(2616,748)(2651,749)
	(2687,749)(2715,750)(2743,752)
	(2773,753)(2804,755)(2835,757)
	(2868,759)(2903,761)(2939,764)
	(2977,767)(3017,771)(3059,775)
	(3103,779)(3149,783)(3198,788)
	(3249,794)(3302,799)(3358,805)
	(3416,811)(3475,818)(3537,825)
	(3600,832)(3663,839)(3727,846)
	(3791,854)(3853,861)(3914,868)
	(3971,875)(4026,881)(4076,887)
	(4121,892)(4161,897)(4195,901)
	(4223,904)(4245,907)(4262,909)
	(4274,910)(4281,911)(4285,912)(4287,912)
\path(162,12)(164,12)(168,13)
	(175,14)(187,15)(204,17)
	(226,20)(254,23)(288,27)
	(328,32)(373,37)(423,43)
	(478,49)(535,56)(596,63)
	(658,70)(722,78)(786,85)
	(849,92)(912,99)(974,106)
	(1033,113)(1091,119)(1147,125)
	(1200,130)(1251,136)(1300,141)
	(1346,145)(1390,149)(1432,153)
	(1472,157)(1510,160)(1546,163)
	(1581,165)(1614,167)(1645,169)
	(1676,171)(1706,172)(1734,174)
	(1762,174)(1798,175)(1833,176)
	(1868,176)(1902,176)(1935,176)
	(1968,175)(2001,175)(2034,174)
	(2066,172)(2099,171)(2132,169)
	(2164,168)(2197,166)(2229,164)
	(2262,162)(2295,160)(2327,158)
	(2360,156)(2392,155)(2425,153)
	(2458,152)(2490,150)(2523,149)
	(2556,149)(2589,148)(2622,148)
	(2656,148)(2691,148)(2726,149)
	(2762,149)(2790,150)(2818,152)
	(2848,153)(2879,155)(2910,157)
	(2943,159)(2978,161)(3014,164)
	(3052,167)(3092,171)(3134,175)
	(3178,179)(3224,183)(3273,188)
	(3324,194)(3377,199)(3433,205)
	(3491,211)(3550,218)(3612,225)
	(3675,232)(3738,239)(3802,246)
	(3866,254)(3928,261)(3989,268)
	(4046,275)(4101,281)(4151,287)
	(4196,292)(4236,297)(4270,301)
	(4298,304)(4320,307)(4337,309)
	(4349,310)(4356,311)(4360,312)(4362,312)
\path(1962,1737)(1960,1735)(1957,1732)
	(1950,1725)(1940,1715)(1928,1701)
	(1912,1685)(1895,1666)(1876,1645)
	(1857,1622)(1837,1599)(1818,1574)
	(1800,1549)(1784,1524)(1769,1498)
	(1756,1472)(1746,1445)(1739,1417)
	(1736,1389)(1737,1362)(1743,1337)
	(1754,1315)(1767,1296)(1784,1281)
	(1802,1267)(1822,1256)(1843,1247)
	(1866,1240)(1889,1234)(1913,1228)
	(1936,1224)(1958,1221)(1979,1218)
	(1997,1216)(2012,1214)(2023,1213)
	(2031,1212)(2035,1212)(2037,1212)
\path(1962,1362)(1960,1360)(1957,1357)
	(1951,1351)(1942,1341)(1931,1329)
	(1918,1315)(1904,1298)(1889,1280)
	(1875,1260)(1860,1240)(1847,1219)
	(1835,1196)(1825,1172)(1817,1147)
	(1812,1120)(1810,1091)(1812,1062)
	(1819,1030)(1831,1001)(1846,976)
	(1863,954)(1882,934)(1903,916)
	(1924,901)(1946,886)(1967,873)
	(1987,862)(2005,853)(2019,846)
	(2028,841)(2034,838)(2037,837)
\put(3012,1437){\makebox(0,0)[lb]{{$\cX$}}}
\put(1962,362){\makebox(0,0)[lb]{{$\dar$}}}
\end{picture}
}

If you think about the family of stable maps $\PP^1 \to X$ parametrized by $\PP^1 \smallsetminus \{0\}$ given by the embedding of $\PP^1$ in the corresponding fiber, there simply isn't any stable map from a nodal curve that can be fit over the missing point $\{0\}$! The only reasonable thing to fit in there is the fiber itself, which is an orbifold nodal curve. We call these {\em twisted curves}.

\subsection{Twisted curves}
This is what happens in general: degenerations force us to allow stacky (or twisted) structure at the nodes. Thinking ahead about gluing curves we see that we had better allow these structures at markings as well.

A {\em twisted curve} is a gadget as follows:

$$ \begin{array}{ccc}
\Sigma_i & \subset & \cC \\
                &             & \dar\\
                &             &  C.                     
\end{array}$$
\setlength{\unitlength}{0.00083333in}
\begingroup\makeatletter\ifx\SetFigFont\undefined%
\gdef\SetFigFont#1#2#3#4#5{%
  \reset@font\fontsize{#1}{#2pt}%
  \fontfamily{#3}\fontseries{#4}\fontshape{#5}%
  \selectfont}%
\fi\endgroup%
{\renewcommand{\dashlinestretch}{30}
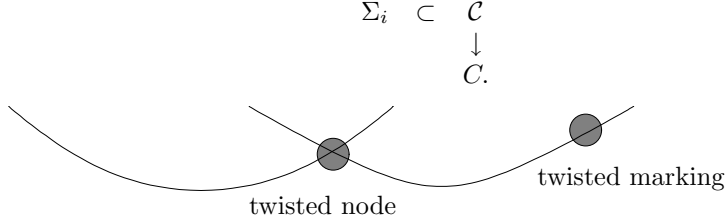
\begin{picture}(4558,715)(0,-10)
\texture{44555555 55aaaaaa aa555555 55aaaaaa aa555555 55aaaaaa aa555555 55aaaaaa 
	aa555555 55aaaaaa aa555555 55aaaaaa aa555555 55aaaaaa aa555555 55aaaaaa 
	aa555555 55aaaaaa aa555555 55aaaaaa aa555555 55aaaaaa aa555555 55aaaaaa 
	aa555555 55aaaaaa aa555555 55aaaaaa aa555555 55aaaaaa aa555555 55aaaaaa }
\put(2037,375){\shade\ellipse{200}{200}}
\put(3612,525){\shade\ellipse{200}{200}}
\path(1512,675)(1514,674)(1517,672)
	(1524,668)(1535,662)(1550,653)
	(1571,641)(1596,627)(1626,610)
	(1661,590)(1700,569)(1743,545)
	(1788,520)(1836,495)(1885,468)
	(1935,442)(1985,416)(2035,391)
	(2084,367)(2132,344)(2178,322)
	(2223,302)(2266,283)(2307,266)
	(2347,251)(2385,237)(2422,224)
	(2457,213)(2491,204)(2525,196)
	(2557,189)(2589,184)(2620,180)
	(2651,177)(2681,176)(2712,175)
	(2743,176)(2773,177)(2804,180)
	(2835,184)(2867,189)(2899,196)
	(2933,204)(2967,213)(3002,224)
	(3039,237)(3077,251)(3117,266)
	(3158,283)(3201,302)(3246,322)
	(3292,344)(3340,367)(3389,391)
	(3439,416)(3489,442)(3539,468)
	(3588,495)(3636,520)(3681,545)
	(3724,569)(3763,590)(3798,610)
	(3828,627)(3853,641)(3874,653)
	(3889,662)(3900,668)(3907,672)
	(3910,674)(3912,675)
\path(12,675)(13,674)(15,673)
	(18,670)(23,666)(30,659)
	(40,651)(53,641)(68,628)
	(86,614)(106,597)(129,579)
	(155,559)(183,538)(212,516)
	(244,493)(277,469)(312,445)
	(348,421)(385,398)(424,374)
	(464,351)(505,329)(547,307)
	(591,286)(637,266)(684,247)
	(734,229)(786,213)(839,198)
	(896,184)(955,173)(1016,163)
	(1080,156)(1146,152)(1212,150)
	(1278,152)(1344,156)(1408,163)
	(1469,173)(1528,184)(1585,198)
	(1638,213)(1690,229)(1740,247)
	(1787,266)(1833,286)(1877,307)
	(1919,329)(1960,351)(2000,374)
	(2039,398)(2076,421)(2112,445)
	(2147,469)(2180,493)(2212,516)
	(2241,538)(2269,559)(2295,579)
	(2318,597)(2338,614)(2356,628)
	(2371,641)(2384,651)(2394,659)
	(2401,666)(2406,670)(2409,673)
	(2411,674)(2412,675)
\put(1512,0){\makebox(0,0)[lb]{{twisted node}}}
\put(3312,150){\makebox(0,0)[lb]{{twisted marking}}}
\end{picture}
}

Here 
\begin{itemize}
\item $C$ is a nodal curve.
\item $\cC$ is a Deligne--Mumford stack with $C$ as its coarse moduli space.
\item Over a node $xy=0$ of $C$, the twisted curve $\cC$ has a chart
 $$[\{uv=0\} / \bmu_r]$$
 where the action of the cyclotomic group $\bmu_r$ is described by 
 $$(u,v) \mapsto (\zeta u , \zeta^{-1} v).$$
 We call this kind of action, with two inverse weights $\zeta, \zeta^{-1}$, a {\em balanced action}. It is necessary for the existence of smoothing of $\cC$! 
In this chart, the map $\cC \to C$ is given by $x = u^r, y = v^r$.
\item At a marking, $\cC$ has a chart $[\AA^1/\mu_r]$, with standard action $u \mapsto \zeta u $, and the map is $x = u^r$.
\item The substack $\Sigma_i$ at the $i$-th marking is locally defined by $u=0$. This stack $\Sigma_i$ is canonically an \'etale gerbe banded by $\bmu_r$.
\end{itemize}
Note that we introduce stacky structure only at isolated points of $C$ and never on whole components. Had we added stack structures along components, we would get in an essential manner a $2$-stack, and I don't really know how to handle these.

As defined, twisted curves form a 2-category, but it is not too hard to show it is equivalent to a category, so we are on safe grounds. 

The automorphism group of a twisted curve is a fascinating object - I'll revisit it later.

This notion of twisted curves was developed in \cite{AV}. As we discovered later, a similar idea appeared in Ekedahl's \cite{Ekedahl}.

\subsection{Twisted stable maps}
\begin{definition}
A twisted stable map consists of 
$$( f: \cC \to \cX, \Sigma_1,\ldots,\Sigma_n),$$ where 
\begin{itemize}
\item $\Sigma_i \subset \cC$ gives a pointed twisted curve.
\item $\cC \stackrel{f}{\to} \cX$  is  a representable morphism.
\item The automorphism group $\Aut_\cX(f, \Sigma_i)$ of $f$ fixing $\Sigma_i$ is finite.
\end{itemize} 
\end{definition}
I need to say something about the last two  {\em stability condition}, necessary for the moduli problem being separated. 

{\em Representability} of $f:\cC \to \cX$ means that for any point $x$ of $\cC$ the associated map
$$\Aut (x) \to \Aut(f(x))$$ on automorphisms is injective. So the orbifold structure on $\cC$ is the ``most economical" possible, in that we do not add unnecessary automorphisms. 

The second condition is in analogy with the usual stable map case, and indeed it can be replaced by conditions on ampleness of a suitable sheaf or number of special points on rational and elliptic components. Most conveniently, it is equivalent to the following schematic condition: the map of course moduli spaces
$$f: C \to X$$ is stable. 

But as I defined things I have not told you what an element of $\Aut_\cX(f, \Sigma_i)$ is! In fact, to make this into a stack I need a category of families of such twisted stable maps. 

\begin{definition} A map   from 
$( f: \cC \to \cX, \Sigma_1,\ldots,\Sigma_n)$ over $S$  to $( f': \cC'\to \cX', \Sigma'_1,\ldots,\Sigma'_n)$ over $S'$ 
 is the following:
$$\xymatrix{
\cC \ar[r]_F\ar[d]\ar@/^1.2pc/[rr]^f & \cC' \ar[r]_{f'}\ar[d]& \cX \\
S \ar[r] & S',
}$$
consisting of 
\begin{itemize}
\item a fiber diagram with morphism $F$ as above, and
\item a 2-isomorphism $\alpha: f \to f'\circ F$. 
\end{itemize}
\end{definition}

Note that the notion of automorphisms is more subtle than the case of  stable maps to a scheme, even if $\cC$ is a scheme. For instance, in the case $\cX = \ocM_g$, a map $\cC \to \cX$ is equivalent to a fibered surface $S \to C$ with fibers of genus $g$, and $S$ can easily have automorphisms acting on the fibers and keeping $C$ fixed, for instance if the fibers are hyperelliptic!

\subsection{Transparency 25: The stack of twisted stable maps}
\subsubsection{The stack of twisted stable maps}
The first nontrivial fact that we have here is the following: the collection of twisted stable maps is again a 2-category, simply because twisted curves are naturally a 2-subcategory of the 2-category of stacks. But there is a simple lemma saying that every 2-morphism between 1-morphisms is unique  and invertible when it exists. This is precisely the condition guaranteeing the following:

\begin{fact}
The 2-category of twisted stable maps is equivalent to a category.
\end{fact}

Here we come to a sticky point: This category is a generalization of $\ocM_{g,n}(X,\beta)$. But in our applications we wish to sometimes insert moduli spaces denoted $\ocM$ in there for $\cX$ - our original application was $\cX  =\ocM_{1,1}$! We made a decision to avoid confusion and denote the category $$\cK_{g,n}(\cX,\beta),$$ after Kontsevich. Some people have objected quite vocally, but I think the choice is sound and the objections are not too convincing and a bit too late for us. But you are perfectly welcome to use notation of your choice.

  The main result is:
  
  \begin{theorem}
The category   $\cK_{g,n}(\cX,\beta)$ is a proper Deligne--Mumford stack with projective coarse moduli space.
  \end{theorem}
  
  Here $\beta$ is simply the class of the curve $f_*[C]$ on the coarse moduli space $X$.

\subsubsection{On the proof}
To prove this, Vistoli and I had to go through several chambers of hell. As far as I can see, the complete symplectic proof (using Fukaya-Ono) is not easier. 

Our main difficulty was the fact that a basic tool like Hilbert schemes was not available for our construction. Now a much better approach, due almost entirely to Martin Olsson, is available. I will sketch is now, because I think it is beautiful. Readers who are not keen on subtleties of constructing moduli stacks might prefer to skip this and take the theorem entirely on faith.
  
  Olsson's proof has the following components:
  
 \begin{itemize}
 \item He first constructs rather explicitly the stack of twisted curves with its universal family \cite{Olsson-log}:
 $$\begin{array}{c} \fC^\tw_{g,n} \\ \dar \\  \fM^\tw_{g,n}. 
 \end{array}$$
 \item In great generality he constructs \cite{Olsson-Hom} a stack of morphisms between two given stacks, and identifies the substack 
 $$\cK_{g,n}(\cX, \beta)  \ \ \subset \ \ \Hom_{\fM^\tw_{g,n}}( \fC^\tw_{g,n}, \cX).$$
 \item He further shows  in the same paper that  when passing to coarse moduli spaces, the natural morphism $$\Hom_{\fM^\tw_{g,n}}( \fC^\tw_{g,n}, \cX) \ \ \lrar \ \ \Hom_{\fM_{g,n}^\tw}( \fC_{g,n}, X) $$ is of finite type, implying the same for $\cK_{g,n}(\cX, \beta) \to \ocM_{g,n}(X, \beta)$.
 \item To prove properness one can use the valuative criterion, whose proof in \cite{AV} is appropriate. In the same paper one counts and sees that $\cK_{g,n}(\cX, \beta) \to \ocM_{g,n}(X, \beta)$ has finite fibers, implying projectivity. 
 \end{itemize} 

\subsection{Twisted curves and roots}

Martin Olsson constructs the stack of twisted curves in general using logarithmic structures \cite{Olsson-log}. This is a very nice construction, but it would have been too much to introduce yet another big theory in these lectures. What I want to do here is describe a variant of this using root stacks, which works nicely in the case of tree-like curves - i.e. where the dual graph is a tree, equivalently every node separates the curve in two connected components.

The construction was first invented by Angelo Vistoli, but his treatment has not yet appeared in print. I lectured on this at ICTP, but did not include in the lecture notes.  Charles Cadman discovered this construction independently and used it to great advantage in his thesis \cite{Cadman}, where a treatment is published.
  
\begin{definition}
Consider a scheme $X$, a line bundle $L$, a section $s \in \Gamma(X,L)$, and a positive integer $r$.
Define a  stack $$\sqrt[r]{(L/X, s)}$$ whose objects over a scheme $Y$ are
$(f: Y \to X, M,\phi, t)$ where  
\begin{itemize} 
\item $M$ is a line bundle on $Y$ and $t\in \Gamma(Y,M)$,
\item $\phi: M^{\otimes r} \stackrel{\sim}{\lrar} L$, and
\item $\phi(t^r) = s$.
\end{itemize} 
\end{definition}
Arrows are fiber diagrams as usual.

For a Cartier divisor $D$, Vistoli uses the notation 
$$\sqrt[r]{(X,D)}:=\sqrt[r]{(\cO_X(D)/X, \boldsymbol{1}_D)}.$$
Cadman uses the notation $X_{D,r}$.


This stack $\sqrt[r]{(X,D)}$ or $X_{D,r}$ is isomorphic to $X$ away from the zero set $D$ of the section, and canonically introduces a stack structure with index $r$ along $D$, which is ``minimal" if $D$ is smooth.  This immediately enables us to define the stacky structure of a twisted curve at a marking starting with the coarse curve: 
$$(C, p) \ \ \rightsquigarrow \ \ \cC = \sqrt[r]{(C,p)} = C_{p,r}. $$

The case of a node is more subtle, and is best treated universally. Here we need to assume that the nodes are separating to use root stacks directly, otherwise one needs either subtle descent or logarithmic structures.

Assume given 
\begin{itemize} 
\item a versal deformation space  of nodal curves $C \to V$, with $V$ a polydisk or a strictly henselian scheme
\item $D \subset V$ the smooth divisor where a particular node in the fibers is preserved,
\item $Z \subset C$ the locus of these nodes, assumed separating, and
\item $E_1, E_2\subset C$ the two connected components of the preimage of $D$ separated by $Z$.  
\end{itemize}
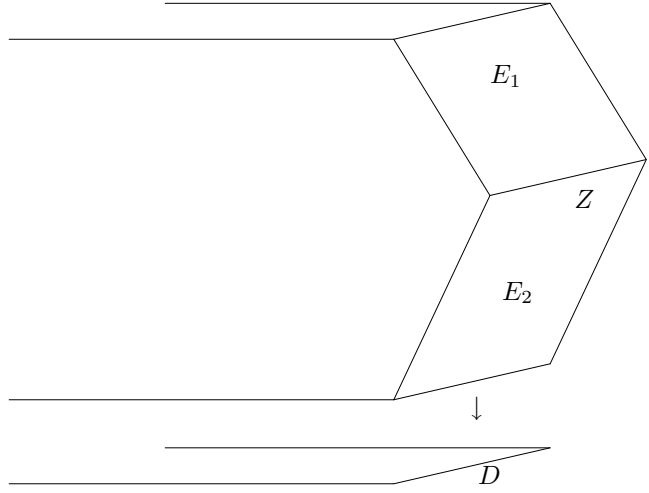
\begin{center}
\setlength{\unitlength}{0.00083333in}
\begingroup\makeatletter\ifx\SetFigFont\undefined%
\gdef\SetFigFont#1#2#3#4#5{%
  \reset@font\fontsize{#1}{#2pt}%
  \fontfamily{#3}\fontseries{#4}\fontshape{#5}%
  \selectfont}%
\fi\endgroup%
{\renewcommand{\dashlinestretch}{30}
\begin{picture}(4005,3039)(0,-10)
\path(12,12)(2412,12)
\path(2412,12)(3387,237)
\path(987,237)(3387,237)
\path(987,3012)(3387,3012)
\path(3387,3012)(3987,2037)(3387,762)
\path(2412,2787)(3012,1812)(2412,537)
\path(12,2787)(2412,2787)
\path(12,537)(2412,537)
\path(2412,2787)(3387,3012)
\path(3012,1812)(3987,2037)
\path(2412,537)(3387,762)
\put(2937,12){\makebox(0,0)[lb]{{$D$}}}
\put(3537,1737){\makebox(0,0)[lb]{{$Z$}}}
\put(3012,2487){\makebox(0,0)[lb]{{$E_1$}}}
\put(3087,1137){\makebox(0,0)[lb]{{$E_2$}}}
\put(2887,412){\makebox(0,0)[lb]{{$\dar$}}}
\end{picture}
}
\end{center}

We have the structure morphism $V \to \fM :=\fM_{g,n}$. Denote by $\fM_r^\tw$ the locus in $\fM_{g,n}^\tw$ where the given node is given stacky structure of index $r$, and  $\fC_r^\tw$ the universal twisted curve. Then we have
\begin{align*}
  V  \ \mathop{\times}\limits_{\fM} \ \fM_r^\tw  \ &= \ \ \sqrt[r]{(V,D)}  \\
 V \ \mathop{\times}\limits_{\fM} \ \fC_r^\tw \ \ &= \ \
 \sqrt[r]{(C,E_1)}\ \mathop{\times}\limits_C \ \sqrt[r]{(C,E_2)}.
 \end{align*}
 
 One can pore over these formulas for a long time to understand them. One thing I like to harvest from the first formula is a description of automorphisms: since  $\fM_r^\tw \to \fM$ is birational, but the versal deformation is branched with index $r$ over $D$, this branching is accounted for by automorphisms of the twisted curve. We deduce that the automorphism group of a twisted curve fixing $C$ is
 $$\Aut_C (\cC)\  = \  \prod_{s\in \Sing \cC} \Gamma_s,$$
 where $\Gamma_s \simeq \bmu_{r_s}$ is the stabilizer of the corresponding node.

These automorphisms acting trivially on $C$ are completely absent from the simple-minded orbifold picture, as in \cite{Chen-Ruan}. Alessio Corti calls them {\em ghost automorphisms}, and their understanding is a key to the paper \cite{Abramovich-Corti-Vistoli}.

\subsection{Valuative criterion for properness}
\subsubsection{Vistoli's Purity Lemma}
The key ingredient   in proving the valuative criterion for properness is the following lemma, stated and proven by Vistoli:

\begin{lemma}[Vistoli]
Consider the following commutative diagram:
$$\xymatrix{
 		& 				& \cX \ar[d]\\ 
{\ \ U\ \ }\ar@{^{(}->}[r]\ar[rru]^{{f}_U} &
		S\ar@{.>}[ru]|{\exists?\, f}
                          \ar[r]_{\overline{f}} &
						X}$$
where 
\begin{itemize}
\item $\cX$ is a Deligne--Mumford stack and $\pi: \cX \to X$ is the coarse moduli space map,
\item $S$ is a smooth variety  and $U \subset S$ is open with compliment of codimension $\geq 2$.
\end{itemize}
Then there exists  $f: S \to \cX$ making the diagram commutative, unique up to a unique isomorphism.
\end{lemma} 

I love this lemma (and its proof). It seems that the lemma loves me back, as it has carried me through half my career!

Related ideas appered in Mochizuki's \cite{Mochizuki}.

\subsubsection{Proof of the lemma}
For the proof, we may assume (working analytically) that $S$ is a $2$-disc and $U$ is a punctured $2$-disc. Hence $U$ is simply connected. We can make $\cX$ smaller as well and assume that $\cX=[V/G]$ where $G$ is a finite group fixing the origin, chosen so that its image in $X$ is the image via $\bar f$ of $S\smallsetminus U=p$.
$$\xymatrix{V\times_XU \ar[r]\ar[d] & V\times_XS \ar[r]\ar[dd] & V \ar[d]\\ 
 \cX\times_XU \ar[d]& 	& \cX \ar[d]\\ 
{\ \ U\ \ }\ar@/_.8pc/[u]_{{f}_U}
\ar@{.>}@/^2.5pc/[uu]^{\exists}
\ar@{^{(}->}[r]
&
		S\ar@{.>}[ru]|{\exists?\, f}
                          \ar[r]_{\overline{f}} &
						X}$$

 Note that $V\times_XU\to \cX\times_XU$ is \'etale and proper, hence the section $U\to \cX\times_XU$ lifts to $V\times_XU$ because $U$ is simply connected. 
Let $\bar U$ be the closure of the image of $U$ in $V\times_XS$. 
$$\xymatrix{ & \bar U\ar[d]\ar[r] & \cX \ar[d]\\ 
{\ \ U\ \ }
\ar@/^1.5pc/[ur]
\ar@{^{(}->}[r]
&
		S\ar@{.>}@/^0.8pc/[u]^{\exists} \ar@{.>}[ru]|{\exists \, f}
                          \ar[r]_{\overline{f}} &
						X}$$

We obtain by projection $\bar U\to S$ which is finite birational and since $S$ is smooth  it has a section. Composing with $\bar U\to \cX$ yields $f$.

\subsubsection{Proof of the valuative criterion}
Let me sketch how the valuative criterion for properness  of $\cK_{g,n}(\cX,\beta)$ follows. Consider a punctured smooth curve $V \subset B$ and a family of stable maps 
$$\xymatrix{
C_V\ar[r]^f\ar[d] & \cX \\V.
}$$ 
By properness    of $\ocM_{g,n}(X,\beta)$ we may assume this extends to a family of stable maps
$$\xymatrix{
C \ar[r]^f\ar[d] & X \\B.
}$$ 
Using properness of $\cX$ and base change (and abhyankar's lemma on fundamental groups - characteristic 0 is essential!) we may assume that the map $C \to X$ lifts to $C \dashrightarrow \cX$ defined on a neighborhood of the generic point of every component of the fiber  $C_0$.  So it is defined on an open set $U\subset C$ whose complement has codimension 2. The purity lemma says that this map extends on the regular locus of $C$. 

 We are left to deal with singular points, where $C$ is described locally as $xy = t^r$ for some $r$. The purity lemma says this extends over the local universal cover given by $uv = t$, where $u^r=x$ and $v^r = y$, and the fundamental group of the punctured neighborhood  is cyclic of order $r$. Of course $r$ changes when you do base change, but there is a smallest $r'$ such that  the map extends over the cover given by $u^{r'} =x, v^{r'} = y$, and this one is representable.  

\section{Gromov--Witten classes}
\subsection{Contractions}
Many of the features of the stacks of stable maps are still true for twisted stable maps. For instance, of one has a morphism $f: \cX \to \cY$ and $m< n$, then as long as either $2g-2+ m > 0$  or $f_*\beta \neq 0$ we have a canonical map
$$\cK_{g,n}(\cX, \beta) \to \cK_{g,m}(\cY, f_\beta), $$
the construction of which is just a bit more subtle than the classical stabilization procedure.

\subsection{Gluing and rigidified inertia}
Much more subtle is the issue of gluing, and the related evaluation maps. To understand it we consider a nodal twisted curve with a separating node:
$$ \cC = \cC_1 \mathop\sqcup\limits^\Sigma \cC_2.$$

As expected, one can prove that $\cC$ is a coproduct in a suitable stack-theoretic sense, and therefore 
$$\Hom(\cC, \cX) \ \ = \ \ \Hom(\cC_1, \cX)\mathop\times\limits_{\Hom(\Sigma, \cX)} \Hom(\cC_2, \cX).$$ 
but $\Sigma$ is no longer a point but a  gerbe! We must ask 

\begin{itemize} 
\item How can we understand $\Hom(\Sigma, \cX)$?
\item What is the universal picture?
\end{itemize}

Since $\Sigma$ is a gerbe banded by $\bmu_r$, the nature of  $\Hom(\Sigma, \cX)$ definitely depends on $r$. 

\begin{definition} Define a category
$$\overline \cI(\cX) \ \ = \ \ \coprod_r \overline \cI_r(\cX), $$
where each component has objects 
$$\overline \cI_r(\cX) (T) \ \ = \ \ \left\{ \begin{array}{ccc} \cG & \stackrel{\phi}{\lrar} & \cX \\ \dar && \\ T && \end{array}\right\}$$
where  
\begin{itemize}
\item $\cG \to T$ is a gerbe banded by $\mu_r$, and
\item $\phi: \cG \to \cX$ is representable.
\end{itemize}
\end{definition}

A priori this is again a 2-category, but again and for a different reason, it is equivalent to a category. In fact we have
\begin{theorem}
$\overline \cI(\cX)$ is a Deligne--Mumford stack.
\end{theorem}
There is a close relationship between $\overline \cI(\cX)$ and the inertia stack $\cI(\cX)$. In fact it is not too difficult to see that there is a diagram 
$$\xymatrix{\cI(\cX) \ar[d]\ar[r] & \cX \\
\overline \cI(\cX)
}$$
making $\cI(\cX)$ the universal gerbe over $\overline \cI(\cX)$!
The stack $\overline \cI_r(\cX)$ can be constructed as the rigidification  of the order $r$ part $$\cI_r(\cX) = \{(\xi,g) | g\in \Aut \xi, g \mbox{ of order } r\} $$ of the inertia stack by removing the action of the cyclic group $\<g\>$ of order $r$ from the picture. This is analogous to the construction of Picard scheme, where the $\CC^*$ automorphisms are removed from the picture. This is discussed in \cite{Abramovich-Corti-Vistoli} and in \cite{Romagny}, and the notation is
$$\overline \cI_r(\cX) = \cI_r(\cX) \thickslash \bmu_r.$$
In a 2-categorical sense, this is the same as saying that the group-stack $\cB \bmu_r$ acts 2-freely on $\cI_r(\cX)$ and in fact
$$\overline \cI_r(\cX) = \cI_r(\cX) / \cB \bmu_r,$$
but this tends to make me dizzy.

The name we give $\overline \cI(\cX)$ is {\em the rigidified inertia stack}.

The following example may be illuminating: consider the global quotient stack $\cX = [Y/G]$. Then a simple analysis shows:
$$ \cI(\cX) \ \ = \ \ \coprod_{(g)}\left[\ Y^g\ \Big/ \ C(g)\ \right],
$$
where the union is over conjugacy classes $(g)$ and $C(g)$ denotes the centralizer of $g$. But by definition the cyclic group $\<g\>$ acts trivially on $Y^g$, therefore the action of $C(g)$ factors through $$\overline{C(g)} := C(g) / \<g\>.$$ We have
$$ \overline \cI(\cX) \ \ = \ \ \coprod_{(g)}\left[\ Y^g\ \Big/ \ \overline{C(g)}\ \right].
$$

\subsection{Evaluation maps}
We now have a natural evaluation map
\begin{align*}
\cK_{g,n}(\cX, \beta)  & \stackrel{e_i}{ \lrar}\ \ \  \overline \cI(\cX) \\
(f:\cC \to \cX, \Sigma_1,\ldots \Sigma_n) &\ \mapsto \ \ \ \Sigma_i \stackrel{f_{\Sigma_i}}{\lrar} \cX
\end{align*}

One point needs to be clarified: the gerbe $\Sigma_i$ needs to be banded in order to give an object of $\overline \cI(\cX)$. This is automatic - as $\bmu_r$ canonically acts on the tangent space of $\cC$ at $\Sigma_i$, this gives a canonical identification of the automorphism  group of a point on $\Sigma_i$ with $\bmu_r$!

There is a subtle feature we have to add here: when gluing curves we need to make the glued curve balanced. Therefore on one branch the banding has to be inverse to the other. We wire this into the definitions as follows. There is a natural involution
\begin{align*}
 \cI(\cX) & \stackrel{\iota}{\lrar}  \cI(\cX)\\
 (\xi, \sigma) & \mapsto (\xi, \sigma^{-1})
\end{align*}
which gives rise to an involution, denoted by the same symbol,
$\iota:  \overline \cI(\cX) \to  \overline \cI(\cX)$. 

On the level of gerbes, it sends $\cG \to \cX$ to itself, but for an object $y$ of $\cG$ changes the isomorphism $\bmu_r \simeq \Aut (Y)$ by composing with the homomorphism  $\bmu_r \to \bmu_r$ sending $\zeta_r \mapsto \zeta_r^{-1}$. 

We define a new {\em twisted evaluation map}
$$
\check e_i = \iota \circ e_i: \cK_{g,n}(\cX, \beta)  \lrar  \overline \cI(\cX). 
$$

\subsection{The boundary of moduli}

We can finally answer the tow questions asked before:

Given a nodal curve $ \cC = \cC_1 \mathop\sqcup\limits^\Sigma \cC_2,$
we have
$$\Hom(\cC, \cX) \ \ = \ \ \Hom(\cC_1, \cX)\mathop\times\limits_{\ \overline \cI(\cX)\ } \Hom(\cC_2, \cX).$$
where the fibered product is with respect to $\check e_\blacktriangleright: \cC_1 \to \overline \cI(\cX)$ on the left and $e_\blacktriangleleft: \cC_2 \to \overline \cI(\cX)$ on the right.

We can apply this principle on universal curves, and obtain, just as in the case of usual stable maps, and obtain a morphism

$$ \cK_{g_1, n_1+\blacktriangleright}(\cX, \beta_1) \mathop\times\limits_{\ \overline \cI(\cX)\ }  \cK_{g_2, n_2+\blacktriangleleft}(\cX, \beta_2) \lrar \cK_{g, n}(X, \beta), $$ with the fibered product over the twisted evaluation map $\check e_{\blacktriangleright}$ on the left and  the non-twisted $e_{\blacktriangleleft}$ on the right. Again this is compatible with the other evaluations, e.g. for $i\leq n_1$ we have a commutative diagram 

$$\xymatrix{ \cK_{g_1, n_1+\blacktriangleright}(\cX, \beta_1) \mathop{\times}\limits_{\ \overline \cI(\cX)\ }  \cK_{g_2, n_2+\blacktriangleleft}(\cX, \beta_2) \ar[r]\ar_{\pi_1}[d] &\cK_{g, n}(\cX, \beta)\ar^{e_i}[d]\\
\cK_{g_1, n_1+\blacktriangleright}(\cX, \beta_1)\ar^{e_i}[r]& {\ \overline \cI(\cX)\ }.
} $$

We now come to a central observation of orbifold Gromov--Witten theory:
\begin{quote} \bf
Since evaluation maps lie in $\overline \cI(\cX)$, Gromov--Witten classes operate on the cohomology of $\overline \cI(\cX)$, and not of $\cX$!!
\end{quote}

Of course $\cX$ is part of the picture: $$\cX = \overline \cI_1(\cX) \subset \overline \cI(\cX).$$ The other pieces of $\overline \cI(\cX)$ are known as {\em twisted sectors} (yet another meaning of ``twisting" in geometry), and arise in various places in mathematics and physics for different reasons. But from the point of view of Gromov--Witten theory, the reason they arise is just the observation above, which comes down to the fact that a nodal twisted curve is glued along a gerbe, not a point.

\subsection{Orbifold Gromov--Witten classes}
First, an easy technicality: we have a locally constant function 
\begin{align*}
r:  \overline \cI(\cX) & \lrar \ZZ \\ 
\intertext{defined by sending }
\overline \cI_r(\cX) & \mapsto r
\end{align*}
On the level of the non-rigidified inertia stack $\cI(\cX)$, it simply means that
$$ (\xi,\sigma) \ \ \mapsto\ \ \mbox{ the order of  } \sigma.$$

A locally constant function gives an element of $H^0( \cI(\cX), \ZZ)$, and therefore we can multiply any cohomology class by this $r$.

Now, as observed above, Gromov--Witten theory operates on $H^*(\overline \cI(\cX))$. 
Because of its role in Gromov--Witten theory, the cohomology space $H^*(\overline \cI(\cX))$ got a special name - it is commonly known as the orbifold cohomology of $\cX$, denoted by 
$$H^*_{orb}(\cX) := H^*(\overline \cI(\cX)).$$

We again simplify notation: $\cK = \cK_{0,n+1}(\cX, \beta)$, and take $\gamma_i \in H^*_{orb}(\cX)^{even}$ to avoid sign issues.

We now define {\em Gromov--Witten classes}:
\begin{definition}
$$\<\gamma_1,\ldots,\gamma_n,*\>^\cX_\beta\ :=\ \  r\cdot \check e_{n+1\,*}(e^*(\gamma_1\cup\ldots\cup\gamma_n)\cap[\cK]^\vir))\ \ \in\ \  H^*_{orb}(\cX).$$
\end{definition}

I'm going to claim that WDVV works for these classes just as before, but any reasonable person will object - this factor of $r$ must enter somewhere, and I had better explain why it is there and how it works out!

There are two ways I want to answer this. First - on a formal level, this factor $r$ is needed because of the following: consider the stack  $\fD(12|34)^\tw\to \fM_{0,4}^\tw$ consisting of twisted nodal curves where markings numbered $1,2$  are separated from $3,4$ by a node. This
is not the fibered product  $(\fM_{0,3}\times \fM_{0,3})\times_{\fM_{0,4}}  \fM_{0,4}^\tw$, and there is exactly a multiplicity $r$ involved - I explain this below. (In \cite{AGV0} we used a slightly different formalism, replacing the moduli stacks $\cK_{g,n}(\cX,\beta)$ with the universal gerbes, and a different correction is necessary).

The second answer is important at least from a practical point of view. The Gromov--Witten theory of $\cX$ wants to behave as if there is a lifted evaluation map $\tilde e_i: \cK_{g,n}(\cX,\beta) \to \cI(\cX)$, lifting 
$e_i: \cK_{g,n}(\cX,\beta) \to \ocI(\cX)$. 
$$\xymatrix{ & \cI(\cX)\ar[d]^{\pi}\\
\cK_{g,n}(\cX,\beta)\ar[r]_{e_i}\ar@{.>}[ru]|{\not\exists\, \tilde e_i} &  \ocI(\cX)
}$$
Of course pulling back gives a multiplicative  isomorphism $$\pi^*: H^*(\ocI(\cX),\QQ) \to H^*(\cI(\cX),\QQ),$$
but a cohomological lifting of, say $\tilde e_{i\,*}$ of $e_{i\,*}$ is obtained rather by composing  with the non-multiplicative  isomorphism 
$$(\pi_*)^{-1}: H^*(\ocI(\cX),\QQ) \to H^*(\cI(\cX),\QQ),$$ and $(\pi_*)^{-1} = r \cdot \pi^*$, 
so $$\tilde e_{i\,*} = (\pi_*)^{-1} \circ e_{i\,*} = r \pi^*\circ e_{i\,*}.$$
Similarly we define  $$\tilde e_i^* = e_i^*\circ (\pi^*)^{-1}.$$ Since $(\pi^*)^{-1}= r\cdot \pi_*$ we can also write $\tilde e_i^* = r\cdot e_i^*\circ \pi_*.$\footnote{Notice that the latter formula was written incorrectly in \cite{AGV0}, section 4.5. We are indebted to Charles Cadman for noting this error. He also used these ``cohomological evaluation maps" very effectively in his work.} The factor $r$ is cancelled out beautifully in WDVV, and one has all the formula as in the classical case. I'll indicate how this works in an example below.

I must admit that I am not entirely satisfied with this situation. I wish there were some true map standing for $\tilde e_i$. This might be very simple (why not replace $\riX$ by $\riX \times \cB\bmu_r$?), but the 2-categorical issues make my head spin when I think about it.

I also need to say something $[\cK]^\vir$,.

\subsection{Fundamental classes}
Let us see now where the differences and similarities are in the statement and formal proof. 

The virtual fundamental class is again defined by the relative obstruction theory  $$E = \bR\pi_*f^*T_{\cX}$$ coming from the diagram of the universal curve
 $$\xymatrix{ \cC \ar_\pi[d] \ar[r]^f & \cX \\ \cK.
 }$$

It satisfies the Behrend-Fantechi relationship

$$\sum_{\beta_1+\beta_2=\beta}\sum_{A\sqcup B = I}
\Delta^!\left([\cK_1]^\vir\times [\cK_2]^\vir\right) \ \ = \ \ \mathfrak{gl}^![\cK]^\vir. $$

as in the diagram

$$\xymatrix{
\cK_1 \times \cK_2 \ar[d] & \cK_1 \times_{\riX} \cK_2 \ar[l]\ar[d]\ar@{^(->}[r] & \sqcup \cK_1\times_{\riX} \cK_2\ar[r]\ar[d] & \cK\ar[d]\\
\riX \times \riX & \riX \ar^{\Delta}[l] &  \fD^\tw(12|3\bullet) \ar^{\mathfrak{gl}}[r] & \fM_{0,4}.}$$

The proof of this relationship is almost the same as that of Behrend-Fantechi, with one added ingredient: in comparing $E$ with $E_i = \bR\pi_{i\,*}f_i^*T_\cX$ coming from  $$\xymatrix{ \cC_i \ar_\pi[d] \ar[r]^f & \cX \\ \cK_i, 
 }$$ the ``difference" is accounted for by $\pi_{\Sigma\,*}f^*_\Sigma T_\cX$ as in $$\xymatrix{ \Sigma \ar_\pi[d] \ar[r]^f & \cX \\ \cK_1\times_{\riX}\cK_2. 
 }$$
 The crucial step in proving the basic relationship is the following ``tangent bundle lemma":
 
 \begin{lemma}
 Assume  $e: S \to \riX$ corresponds to $$\xymatrix{ \cG \ar_\pi[d] \ar[r]^f & \cX \\ S.
 }$$
 Then there is a canonical isomorphism $$e^* T_{\riX} \stackrel{\sim}\lrar \pi_* f^* T_\cX.$$
 \end{lemma}

I will not prove the lemma  here, but I wish to avow that it is a wonderful proof (which almost fits in the margins). 

\section{WDVV, grading and computations}
\subsection{The formula}
The WDVV formula says again:
\begin{theorem}
\begin{align*}
&\sum_{\beta_1+\beta_2=\beta} \sum_{A\sqcup B = I}
\Big\<\<\gamma_1,\gamma_2,\, \delta_{A_1},\ldots,\delta_{A_k},\, *\>_{\beta_1},\gamma_3,\, 
\delta_{B_1},\ldots,\delta_{B_m},\, * \Big\>_{\beta_2} \\
=&\sum_{\beta_1+\beta_2=\beta}\sum_{A\sqcup B = I}
\Big\<\<\gamma_1,\gamma_3,\, \delta_{A_1},\ldots,\delta_{A_k},\, *\>_{\beta_1},\gamma_2,\, 
\delta_{B_1},\ldots,\delta_{B_m},\, * \Big\>_{\beta_2}
\end{align*}
\end{theorem}
 
 The fibered product diagram looks like this:

$$\xymatrix{
\cK_1 \times_{\ocI(\cX)}\cK_2\ar^{p_2}[r]\ar_{p_1}[d] & \cK_2\ar^{\check e_\bullet}[r] \ar^{e_\blacktriangleleft}[d] & \ocI(\cX) \\
\cK_1 \ar^{\check e_\blacktriangleright}[r] & \ocI(\cX).
}$$

We now have by definition
\begin{align*}\Big\<\<\gamma_1,\gamma_2,\, \delta_{A},\, *\>_{\beta_1},\gamma_3,\, 
\delta_{B},\, * \Big\>_{\beta_2}& = 
r \check e_{\bullet\,*}\Big(e_\blacktriangleleft^*\left(r \check e_{\blacktriangleright\,*}\xi_1\right)\cap \xi_2\Big) \\
\intertext{and the projection formula gives}
& = (r \check e_\bullet\circ p_2)_*\Big(r p_1^*\eta_1 \cup p_2^* \eta_2\cap \Delta^![\cK_1\times \cK_2]\Big)
\end{align*}
where the $r$ inside the parentheses is pulled back from the glued markings. The overall factor $r$ is the same on all terms of the equation so it can be crossed out in the proof. 

The real difference comes in the divisor diagram, which becomes the following fiber diagram:
 
 $$\xymatrix{
{}\save+<0pt,-0pt>
*+{\coprod^{\vphantom{\riX}} \cK_1 \mathop\times\limits_{\riX}\cK_2}
\ar^{\ell}[rrr] \ar_\phi[d]
\restore
   & &  & \cK\ar[d]\\
\fD(12|3\bullet)^\tw \ar^{i}[r] &\star\ar[d]\ar[rr] && \fM^\tw_{0,4}\ar[d]
 \\
& \fM_{0,3}\times \fM_{0,3} \ar^(.7){j}[r] & \star \ar[r]\ar[d] & \fM_{0,4} \ar[d]
  \\ 
  && \{pt\}\ar[r] &     \ocM_{0,4}.
}$$

Now the map $j$ is still birational, but $i$ has degree $1/r$, which exactly cancels the factor $r$ we introduced in the definition of Gromov--Witten classes.

\subsection{Quantum cohomology and its grading}
\subsubsection{Small quantum cohomology} In order to make explicit calculations with Gromov--Witten classes it is useful to have a guiding formalism. Consider, for instance the {\em small quantum cohomology product}. Let $N^+(X)$ be the monoid  of homology classes of effective curves on $X$. Consider the formal monoid-algebra $\QQ[[N^+(X)]]$ in which we represent a generator corresponding to $\beta\in N^+(X)$ multiplicatively as $\q^\beta$. Define a bilinear map on $QH^*(\cX) := H^*(\riX, \QQ)\otimes \QQ[[N^+(X)]]$ by the following rule on generators $\gamma_i \in H^*(\riX, \QQ)$:
$$\gamma_1*\gamma_2 \ \ = \ \ \sum_{\beta\in N^+(X)}\<\gamma_1,\gamma_2,*\>_\beta \cdot \q^\beta.$$
 This is a skew-commutative ring and as in the non-orbifold case,  WDVV gives its associativity. 

A somewhat simpler ring, but still exotic, is obtained by setting all $\q^\beta=0$ for nonzero $\beta$. The product becomes 
$$\gamma_1\cdot\gamma_2 \ \ = \ \ \<\gamma_1,\gamma_2,*\>_0.$$
 This is the so called {\em orbifold cohomology ring} (or {\em string cohomology}, according to Kontsevich).
We denote it $H^*_{orb}(\cX)$. The underlying group is $H^*(\riX, \QQ)$.

But one very very useful fact is that these  are graded rings, and the grading, discovered by physicists, is fascinating.

\subsubsection{Age of a representation}
A homomorphism $\rho: \bmu_r \to \GG_m$ is determined by an integer $0\leq k \leq r-1$, as $\rho(\zeta_r) = \zeta_r^a$. We define  $$age(\rho):= k/r.$$ 

This extends by linearity to a function on the representation ring $age:R\bmu \to \QQ$. 
\subsubsection{Age of a gerbe in $\cX$}
Now consider a morphism $e:T \to \riX$ corresponding to $$\xymatrix{ \cG \ar^f[r]\ar[d] & \cX \\T.}$$ The pullback $f^*T_\cX$ is a representation of $\bmu_r$, and we can define $age(e) : = age (T_\cX)$. 

We can thus define a locally constant function $age: \riX \to \QQ$, the value on any component being the age of any object $e$ evaluating in this component. 

\subsubsection{Orbifold Riemann-Roch} 
There are numerous justifications for the definition of age in the literature, some look a bit like voodoo. In fact,
the true natural justification for the definition of age is Riemann-Roch on orbifold curves.

Consider a twisted curve $\cC$ and a locally free sheaf $\cE$  on $\cC$.
There is a well defined notion of degree $\deg \cE \in \QQ$. It can be defined by the property that if $D$ is a curve and $\phi:D \to \cC$ is surjective of degree $k$ on each component, then $\deg_{\cC} \cE = \deg_D\phi^*\cE / k$

Riemann-Roch says the following:

\begin{theorem}
$$\chi(\cC,\cE) \ \ =\ \  \operatorname{rank} (\cE)\cdot \chi(\cC,\cO_\cC)\ +\ \deg \cE \ \ - \sum_{p_i\mbox{ marked}} age_{p_i} \cE$$
\end{theorem}
This is easy to show for a smooth twisted curve, and not very difficult for a nodal one as well.

\subsection{Grading the rings}
We now define 
$$H^i_{orb}(\cX) \ \ = \ \ \bigoplus_\Omega H^{i\ -\ 2\,age(\Omega)}(\Omega, \QQ),$$
where the sum is over connected components $\Omega\subset \riX$.

We also define $$\deg \q^\beta \ \ = \ \  2 \ \beta'\cdot c_1(T_\cX)$$
where $\beta'$ is any class satisfying $\pi_*\beta' = \beta$. 

Here is the result:

\begin{theorem}
\begin{enumerate}
\item $H^*_{orb}(\cX)$ is a graded ring.
\item The product on $QH^*(\cX)$ is homogeneous, so it is a ``pro-graded" ring.
\end{enumerate}
\end{theorem}

\subsection{Examples}
\begin{example}
Consider $\cX = \cB G$. In this case $$\iX\  =\  [G/G]\  =\  \coprod_{(g)} \cB(C(g)) $$ where $G$ acts on itself by conjugation. This case is simple as there are no nonzero curve classes  and $QH^*(\cX)  = H^*_{orb}(\cX) = \oplus_{(g)}\QQ $. Since the tangent bundle is zero the ages are all 0, and all obstruction bundles vanish - i.e. the virtual fundamental class equals the fundamental class. Also, since $\iX \to \riX$ has a section, we have  lifted evaluation maps $\tilde e_i:\cK_{0,3}(\cB G, 0 ) \to \iX$.

Let us identify the components of $\cK_{0,3}(\cB G, 0)$: 

There is one  component $\cK_{(g,h)}$ for each conjugacy class of  triple $g,h, (gh)^{-1}$, describing the monodromy of a $G$-covering of $\PP^1$ branched at $0,1,\infty$. The automorphism group of such a cover is $C(g)\cap C(h)$, thus $\cK_{(g,h)}\simeq \cB(C(g)\cap C(h))$.  The third (twisted, lifted) evaluation map is
$$\check {\tilde e}_3 : \cK_{g,h} \to \cB (C(gh))$$ of degree
 $$ \deg \check {\tilde e}_3 = \frac{|C(gh)|}{|C(g)\cap C(h))|}.$$ Therefore we get 
 $$x_{(g)}\cdot x_{(h)}=\sum_{(g,h)}\frac{|C(gh)|}{|C(g)\cap C(h)|}x_{(gh)}$$
where the sum runs over all simultaneous conjugacy classes of the pair $(g,h)$.
In other words, $QH^*_{orb}(BG)=Z(\QQ[G])$, the center of the group ring of $G$.
\end{example}

\begin{example} Maybe the simplest nontrivial case is the following, but already here we can see some of the subtleties of the subject. I'll give these in full detail here - once you get the hang of it it gets pretty fast!

 We take $\cX=\PP(p,1)$ where $p$ is a prime. It is known as the {\em teardrop} projective line, a name inspired by imagining this to be a Riemann sphere with a ``slightly pinched" north pole, with conformal angle $2\pi/p$.
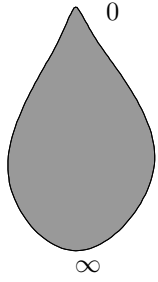
 \begin{center}
 \setlength{\unitlength}{0.00083333in}
\begingroup\makeatletter\ifx\SetFigFont\undefined%
\gdef\SetFigFont#1#2#3#4#5{%
  \reset@font\fontsize{#1}{#2pt}%
  \fontfamily{#3}\fontseries{#4}\fontshape{#5}%
  \selectfont}%
\fi\endgroup%
{\renewcommand{\dashlinestretch}{30}
\begin{picture}(940,1779)(0,-10)
\texture{55888888 88555555 5522a222 a2555555 55888888 88555555 552a2a2a 2a555555 
	55888888 88555555 55a222a2 22555555 55888888 88555555 552a2a2a 2a555555 
	55888888 88555555 5522a222 a2555555 55888888 88555555 552a2a2a 2a555555 
	55888888 88555555 55a222a2 22555555 55888888 88555555 552a2a2a 2a555555 }
\shade\path(610,1436)(623,1418)(637,1399)
	(651,1379)(666,1358)(682,1336)
	(698,1313)(715,1289)(732,1265)
	(749,1239)(766,1212)(783,1185)
	(800,1157)(817,1129)(832,1100)
	(847,1071)(861,1041)(874,1012)
	(886,982)(897,953)(906,924)
	(913,895)(919,866)(924,838)
	(927,810)(928,782)(928,754)
	(926,728)(922,702)(917,675)
	(911,648)(903,621)(893,594)
	(882,566)(870,539)(856,512)
	(841,484)(824,458)(807,431)
	(788,406)(768,381)(748,357)
	(727,335)(705,314)(683,294)
	(660,275)(637,259)(614,244)
	(592,231)(569,219)(546,209)
	(524,202)(501,196)(479,191)
	(458,189)(436,188)(414,190)
	(392,193)(370,198)(348,204)
	(326,213)(304,223)(283,235)
	(261,249)(240,264)(219,281)
	(198,299)(178,319)(159,340)
	(141,363)(124,386)(107,410)
	(92,434)(78,460)(66,485)
	(54,511)(44,537)(35,564)
	(28,590)(22,616)(18,642)
	(14,668)(13,694)(12,722)
	(13,750)(15,779)(19,808)
	(24,838)(31,868)(39,898)
	(48,929)(58,960)(69,992)
	(82,1023)(95,1055)(108,1086)
	(123,1117)(137,1147)(152,1177)
	(167,1205)(182,1233)(197,1261)
	(211,1287)(225,1312)(239,1335)
	(252,1358)(265,1380)(276,1401)
	(288,1421)(304,1453)(320,1483)
	(334,1512)(347,1539)(359,1565)
	(370,1590)(380,1612)(389,1633)
	(397,1651)(405,1667)(411,1680)
	(417,1691)(422,1699)(426,1705)
	(431,1708)(435,1709)(439,1708)
	(443,1706)(447,1701)(452,1695)
	(458,1687)(464,1677)(471,1665)
	(479,1651)(488,1635)(498,1618)
	(509,1599)(520,1578)(533,1557)
	(546,1534)(561,1511)(576,1487)
	(592,1462)(610,1436)
\path(610,1436)(623,1418)(637,1399)
	(651,1379)(666,1358)(682,1336)
	(698,1313)(715,1289)(732,1265)
	(749,1239)(766,1212)(783,1185)
	(800,1157)(817,1129)(832,1100)
	(847,1071)(861,1041)(874,1012)
	(886,982)(897,953)(906,924)
	(913,895)(919,866)(924,838)
	(927,810)(928,782)(928,754)
	(926,728)(922,702)(917,675)
	(911,648)(903,621)(893,594)
	(882,566)(870,539)(856,512)
	(841,484)(824,458)(807,431)
	(788,406)(768,381)(748,357)
	(727,335)(705,314)(683,294)
	(660,275)(637,259)(614,244)
	(592,231)(569,219)(546,209)
	(524,202)(501,196)(479,191)
	(458,189)(436,188)(414,190)
	(392,193)(370,198)(348,204)
	(326,213)(304,223)(283,235)
	(261,249)(240,264)(219,281)
	(198,299)(178,319)(159,340)
	(141,363)(124,386)(107,410)
	(92,434)(78,460)(66,485)
	(54,511)(44,537)(35,564)
	(28,590)(22,616)(18,642)
	(14,668)(13,694)(12,722)
	(13,750)(15,779)(19,808)
	(24,838)(31,868)(39,898)
	(48,929)(58,960)(69,992)
	(82,1023)(95,1055)(108,1086)
	(123,1117)(137,1147)(152,1177)
	(167,1205)(182,1233)(197,1261)
	(211,1287)(225,1312)(239,1335)
	(252,1358)(265,1380)(276,1401)
	(288,1421)(304,1453)(320,1483)
	(334,1512)(347,1539)(359,1565)
	(370,1590)(380,1612)(389,1633)
	(397,1651)(405,1667)(411,1680)
	(417,1691)(422,1699)(426,1705)
	(431,1708)(435,1709)(439,1708)
	(443,1706)(447,1701)(452,1695)
	(458,1687)(464,1677)(471,1665)
	(479,1651)(488,1635)(498,1618)
	(509,1599)(520,1578)(533,1557)
	(546,1534)(561,1511)(576,1487)
	(592,1462)(610,1436)
\put(430,54){\makebox(0,0)[lb]{{$\infty$}}}
\put(625,1629){\makebox(0,0)[lb]{{0}}}
\end{picture}
}
 \end{center}
  It is defined to be the quotient of $\CC^2\smallsetminus\{0\}$ by the $\CC^*$ action $\lambda(x,y)=(\lambda^px,\lambda y)$. In this case $H_2(X)$ is cyclic with a unique positive generator $\beta_1$, and writing $\q = \q^{\beta_1}$ we have $\deg \q=2(p+1)/p$ (in general for the weighted projective space $\PP(a,b)$, with $a, b$ coprime  we have $\deg \q=2(1/a+1/b)$).

We can describe the inertia stack as follows: 
$$\cI(\cX) = \cX \sqcup \coprod_{i=1}^{p-1} \cB (\ZZ/p\ZZ).$$ The rigidified inertia stack is therefore
 $$\ocI(\cX) = \cX \sqcup \coprod_{i=1}^{p-1} \Omega_i,$$ where the rigidified twisted sectors $\Omega_i$ are just points. Since $\Omega_i$ corresponds to the element $i\in \ZZ/p\ZZ$, and $i$ acts on the tangent space  of $\cX$ via $\zeta_p^i$, its age is $i/p$. 
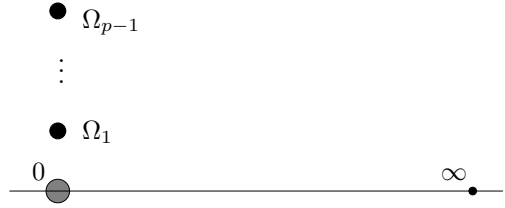
 \begin{center}
\setlength{\unitlength}{0.00083333in}
\begingroup\makeatletter\ifx\SetFigFont\undefined%
\gdef\SetFigFont#1#2#3#4#5{%
  \reset@font\fontsize{#1}{#2pt}%
  \fontfamily{#3}\fontseries{#4}\fontshape{#5}%
  \selectfont}%
\fi\endgroup%
{\renewcommand{\dashlinestretch}{30}
\begin{picture}(3099,1306)(0,-10)
\texture{44555555 55aaaaaa aa555555 55aaaaaa aa555555 55aaaaaa aa555555 55aaaaaa 
	aa555555 55aaaaaa aa555555 55aaaaaa aa555555 55aaaaaa aa555555 55aaaaaa 
	aa555555 55aaaaaa aa555555 55aaaaaa aa555555 55aaaaaa aa555555 55aaaaaa 
	aa555555 55aaaaaa aa555555 55aaaaaa aa555555 55aaaaaa aa555555 55aaaaaa }
\put(312,83){\shade\ellipse{150}{150}}
\put(312,458){\blacken\ellipse{100}{100}}
\put(312,1208){\blacken\ellipse{100}{100}}
\put(2900,83){\blacken\ellipse{50}{50}}
\path(12,83)(3087,83)
\put(462,383){\makebox(0,0)[lb]{{$\Omega_1$}}}
\put(462,1058){\makebox(0,0)[lb]{{$\Omega_{p-1}$}}}
\put(312,758){\makebox(0,0)[lb]{{$\vdots$}}}
\put(150,150){\makebox(0,0)[lb]{{0}}}
\put(2700,150){\makebox(0,0)[lb]{{$\infty$}}}
\end{picture}
}\end{center}
Choose as generator of $H^2(\cX)$ the class $x=c_1(\cO_{\cX}(1))$; we have $\deg x=1/p$. Choose also generators $A_i$ to be each a generator of $H^0(\Omega_i)$ for $1\le i\le p-1$. Since the age is $i/p$, it is positioned in degree $2i/p$ in $H^*_{orb}(\cX)$. We conclude that the orbifold cohomology group, with grading, looks as follows:

$$\begin{array}{rccccccccc}
\text{degree} & 0    &             & 2/p       &          & \ldots &        & 2(p-1)/p && 2 \\
& \QQ& \oplus &\QQ A_1&\oplus&\ldots  & \oplus &\QQ A_{p-1}& \oplus &\QQ x
\end{array}$$

This is a good case to work out the difference between working with $e_i: \cK_{0,3}(\cX, \beta) \to \riX$ and the lifted $\tilde e_i:\cK_{0,3}(\cX, \beta) \to \iX$ (which still exists on  the relevant components of  $\cK_{0,3}(\cX, \beta)$ in this example).

Consider for instance the product $A_i * A_1$ when $i<p-1$. In order to match the degrees, the  result must be  a multiple of $A_{i+1}$, and the only curve class $\beta$ possible is $0$. The component of $\cK_{0,3}(\cX, 0)$ evaluating at $i, 1$ and the inverse $p-i-1$ of $i+1$ is, just as in the first example, the classifying stack of the joint  centralizer of 1 and $i$, which is just $\cB(\ZZ/p\ZZ)$. Since $A_1$ and $A-i$ are fundamental classes, the pullback $e_1^* A_i \cup e_2^* A_1$ is the fundamental class. 

But how about the virtual fundamental class? Here  we use a fundamental fact about the virtual fundamental class: if it a multiple of the fundamental class, then it coincides with the fundamental class.
We calculate $$\check e_{3\,*}(e_1^* A_i \cup e_2^* A_1) = 1/p A_{i+1},$$ since the degree of $\check e_3$ is $1/p$. Since $r=p$ we get 
$$A_i * A_1 = A_i \cdot A_1 = A_{i+1}.$$

Note that, had we used $\tilde e_i$ instead of $e_i$ and replaced $A_i$ by $\tilde A_i = \pi^*A_i$, the change in the degrees of the maps $e_i$ would cancel out with the factor $r$ and we would get the same formula:
$$\tilde A_i * \tilde A_1 = \tilde A_i \cdot \tilde A_1 = \tilde A_{i+1}.$$

Back to the calculation. By associativity, it follows that $A_i*A_j = A_i\cdot A_j = A_{i+j}$ as long as $i+j<p$. 

The next product to calculate is $A_1 *A_{p-1}$. The result must  be a multiple of $x$, again the only class $\beta$ involved is $\beta=0$. Again the virtual fundamental class  is the fundamental class, and we have 
 $$\check e_{3\,*}(e_1^* A_1 \cup e_2^* A_{p-1}) = x,$$ as it has degree $1/p$. This time $r=1$ so
$$A_1 * A_{p-1} = A_1 \cdot A_{p-1} = x.$$

The interesting product is $x*A_1$. It is easy to see that $x\cdot A_1=0$ because there is nothing in degree $2+2/p$ in $H^*_{orb}(\cX)$, and the only possible contribution to $x*A_1$ comes with $\beta=\beta_1$ of degree $2(1+ 1/p)$. What do we see in $\cK_{0,3}(\cX, \beta_1)$? Maps evaluating in $\Omega_1$,  parametrize twisted curves generically of the form $\cX$,  with one (in our case, the second) marking over the orbifold point and two freely roaming around $\cX$.   If we represent $x$ by $1/p \cdot [z]$, with $z$ one of the non-stacky points, then the class $e_1^* x \cup e_1^* A_1$ is represented by a family of twisted stable maps parametrized by the position of the third marking - that is a copy of $\cX$!  The restriction of the virtual fundamental class {\em to the locus $e_1^{-1} z$} is again the fundamental class, and $e_3=\check e_3$ is an isomorphism. We get $[z]*A_1 = \q,$ or $x * A_1 = p^{-1} \q$.  
which completely determines the ring as 
$$QH^*(\cX)\  = \ \QQ[[\q]] [A_1]/(pA_1^{p+1} - \q).$$ 

Again, had we used $\tilde e_i$ instead of $e_i$ and replaced $A_i$ by $\tilde A_i = \pi^*A_i$, the change in the degrees of the maps $e_i$ would cancel out with the factor $r$ and we would get the same ring $\QQ[[\q]] [\tilde A_1]/(p\tilde A_1^{p+1} - \q).$ This is Cadman's preferred formalism.
\end{example}

The more general case of $\PP(a,b)$ will be included in \cite{AGV}.

\subsection{Other work} Already in their original papers, Chen and Ruan gave a good number of examples of orbifold cohomology rings. Much work has been done on orbifold cohomology since, and I have no chance of doing justice to all the contributions. I'll mention a few which come to my rather incomplete memory: Fantechi and G\"ottsche \cite{Fantechi-Gottsche} introduced a method for calculating orbifold cohomology of global quotients and figured out key examples; Uribe \cite{Uribe} also studied symmetric powers;  Borisov, Chen and Smith \cite{Borisov-Chen-Smith} described the orbifold cohomology of toric stacks; Jarvis, Kaufmann and Kimura \cite{JKK2} and Goldin, Holm and Knutson \cite{Goldin-Holm-Knutson} following  Chen and Hu \cite{Chen-Hu} described orbifold cohomology directly without recourse to Gromov--Witten theory.

Computations in Gromov--Witten theory of stacks beyond  orbifold cohomology are not as numerous. C. Cadman \cite{Cadman, Cadman2} computed Gromov--Witten invariants of $\sqrt[r]{(\PP^2, C)}$ with $C$ a smooth cubic, and derived the number of rational plane curves of degree $d$ with tangency conditions to the cubic $C$.  H.-H. Tseng \cite{Tseng} generalized the work of Givental and Coates on the mirror predictions to the case of orbifolds, a subject I'll comment on next.

\subsection{Mirror symmetry and the crepant resolution conjecture}
When string theorists  first came up with the idea of mirror symmetry \cite{COGP}, the basic example was that of the quintic threefold $Q$ and its mirror $X$, which happens to have Gorenstein quotient singularities. Mirror symmetry equates period integrals on the moduli space of complex structures on $X$ with Gromov--Witten numbers of $Q$, and vice versa.  String theorists had no problem dealing with orbifolds \cite{DHVW}, \cite{DHVW2}, but they also had the insight of replacing $X$, or the stack $\cX$, by a crepant resolution $f:Y\to X$. (Recall that a birational map of Gorenstein normal varieties is {\em crepant} if $f^* K_X = K_Y$.  The unique person who could possibly come up with such jocular terms - standing for ``zero discrepancies" - is Miles Reid.) Thus, string theorists posited that period integrals on the moduli space of $Y$ equal Gromov--Witten numbers of $Q$, and vice versa.

One direction - relating period integrals of $Y$ with Gromov--Witten invariants of $Q$, was proven, through elaborate direct computations, independently by Givental \cite{Givental}  and Lian-Liu-Yau \cite{Lian-Liu-Yau}. Unfortunately this endeavor is marred by a disgusting and rather unnecessary controversy, which I will avoid altogether. One important point here is that the period integrals on the moduli space complex structures on $Y$ coincide almost trivially with those of $X$, as any deformation of $X$ persists with the same quotient singularities and the crepant resolution deforms along. So there is not a great issue of passing from $X$, or $\cX$, to $Y$.

The other side is different: do we want to compare period integrals for $Q$ with Gromov--Witten invariants of $\cX$, or of $Y$? The first question that comes up is, does there exist a crepant resolution? The answer for threefolds is yes, one way to show it is using Nakajima's $G$-Hilbert scheme - and my favorite approach is that of Bridgeland, King and Reid \cite{BKR}. In higher dimensions crepant resolutions may not exist and one needs to work directly on $\cX$. It was Tseng in his thesis \cite{Tseng} who addressed the mirror prediction computations directly on the orbifold. 

But the question remains, should we work with $\cX$ or $Y$, and which $Y$ for that matter?
String theorists claim that it doesn't matter - Mirror symmetry works for both. We get the following {\em crepant resolution conjecture}:
\begin{conjecture}[Ruan]
The Gromov--Witten theories of $\cX$ and $Y$ are equivalent.
\end{conjecture}

Moreover, it is conjectured that an explicit procedure of a particular type allows passing from $\cX$ to $Y$ and back.

This was Ruan's original impetus for developing orbifold cohomology!

There is a growing body of work on this exciting subject. See Li--Qin--Wang \cite{Li-Qin-Wang} and Bryan--Graber--Pandharipande \cite{Bryan-Graber-Pandharipande}.

\appendix
\section{The legend of String Cohomology: two letters of Maxim Kontsevich to Lev Borisov}
\renewcommand{\thesubsection}{\thesection.\arabic{subsection}}

\subsection{The legend of String Cohomology}
On December 7, 1995 Kontsevich delivered  a history-making lecture at Orsay, titled {\em String Cohomology}. ``String cohomology" is the name Kontsevich chose to give what we know now, after Chen-Ruan, as {\em orbifold cohomology}, and  Kontsevich's lecture notes  described  the orbifold and quantum cohomology of a global quotient orbifold. Twisted sectors, the age grading, and a version of orbifold stable maps for global quotients are all there.   

Kontsevich never did publish his work on string cohomology. I met him in 1999 to discuss my work with Vistoli on twisted stable maps, and he told me a few hints  about the Gromov--Witten theory of orbifolds (including the fact that the cohomology of the inertia stack is the model for quantum cohomology of $\cX$) which I could not appreciate at the time. He was also aware that Chen and Ruan were pursuing the subject so he had no intention to publish. We do, however, have written evidence in the form of two electronic mail letters he sent Lev Borisov in July 1996. These letters are reproduced below - verbatim with the exception of typographics -   with Kontsevich and Borisov's permission.

What's more - Kontsevich never lectured about string cohomology after all! When I said he delivered  a history-making lecture at Orsay, titled {\em String Cohomology}, I did not lie. The way Fran\c{c}ois Loeser tells the story, this seems like an instance of the Heisenberg Uncertainty Principle: Loeser heard at the time Kontsevich had discovered  a complex analogue of $p$-adic integration, and called Kontsevich to ask about it. In response Kontsevich told him to attend his December 7 lecture at Orsay on  String Cohomology. Thus the lecture Kontsevich did deliver was indeed a history-making  lecture on motivic integration.

For years people speculated what is this ``String Cohomology", which they supposed Kontsevich aimed to develop out of Motivic Integration \ldots

Needless to say, Kontsevich never did publish his work on motivic integration either.  A four-page set of notes is available in \cite{Kontsevich-motivic}. Of course that theory is by now fully developed.

\subsection{The archaeological letters}
Following are the two electronic mail  letters of Kontsevich to Borisov from July 1996. I have trimmed the mailer headings, and, following Kontsevich's request, corrected typographical errors. I also added typesetting commands. Otherwise the text is Kontsevich's text verbatim.

The theory sketched here is the Gromov--Witten theory of a global quotient orbifold, a theory treated in great detail in \cite{Jarvis-Kaufmann-Kimura}.
\newcommand{\Curves}{\operatorname{Curves}}
\newcommand{\PreCurves}{\operatorname{Pre-Curves}}

\subsubsection{Letter of 23 July 1996}\hfill

\noindent
Date: Tue, 23 Jul 96 12:39:42 +0200\\
From: maxim@ihes.fr (Maxim Kontsevich)\\
To: lborisov@msri.org\\
Subject: Re: stable maps to orbifolds\\

Dear Lev,

I didn't write yet anywhere the definition of stable maps to orbifolds,
although I am planning to do it.

Here is it:
\begin{enumerate}
\item  For an orbifold $X=Y/G$ ($G$ is a finite group)
we define a new orbifold $X_1$ as the quotient
of $$Y_1:=\{(x,g)|x \in X, g \in G, gx=x\}$$ by the action of $G$:
$$f(x,g):=(fx,fgf^{-1}).$$ $Y_1$ is a manifold consisting of parts of different
dimensions, including $Y$ itself (for $g=1$).
\item String cohomology of $X$, say $HS(X)$
 are defined as usual rational cohomology of $X_1$. I consider it for a
moment only as a $\ZZ/2$-graded vector space (super space).

One can also easily define $HS$ for orbifolds which are not global quotients
of smooth manifolds.
\item Let $Y$ be an almost complex manifold. Denote by $\Curves(Y)$ the
stack of triples $(C,S,\phi)$ where
\begin{itemize}
\item
 $C$ is a nonempty compact complex curve with may be double points.
 \item   $C$ is not necessarily connected.
\item  $S$ is a finite subset of the smooth part of $C$,
\item  $\phi$ is a holomorphic map from $C$ to $Y$.
\end{itemize}
Stability condition is the absence of infinitesimal automorphisms.

The stack $\Curves(Y)$ consists of infinitely many components which are
products of symmetric powers of usual moduli stacks of stable maps.
\item  We define $\PreCurves(X)$ for $X=Y/G$ as the stack quotient of $\Curves(Y)$ by the
obvious action of $G$.\footnote{This is not the standard notion of stack quotient. The next sentence explains what Kontsevich meant - D.A.}
Precisely it means that we consider quadruples
$(C,S,\phi,A)$ where first 3 terms are as above and $A$ is the action
of $G$ on $(C,S,\phi)$. In other words,
we consider curves with $G$ action and marked points and equivariant maps
to $X$.
 We define $\Curves(X)$ as the subset of $\PreCurves(X)$ consisting of things
where $G$ acts FREELY on $C$ minus $\{$singular points and smooth points$\}$.
This is the essential part of the definition\footnote{Kontsevich omitted here the condition that $G$ should not switch the branches at a node and have balanced action. It is automatic on things that can be smoothed - D.A.}.

Lemma: $\Curves(X)$ is open and closed substack in $\PreCurves(X)$.

As usual we can define the virtual fundamental class of each connected
component
of $\Curves(X)$.
\item For a curve from $\Curves(X)$ we can say when it is ``connected":
when the quotient curve $C/G$ is connected.
For ``connected" curves we define their genus as genus of $C/G$,
and the set of marked points as $S/G$.
We define numbered stable map to $X$ as a ``connected" curve with numbered ``marked points".
\end{enumerate}

The evaluation map from the stable curve with numbered marked points
takes value in the auxiliary orbifold $X_1$ (see 1),2)).

Thus we have a lot of symmetric tensors in $HS(X)$, and they all satisfy
axioms from my paper with Manin. The interesting thing is the grading on $HS$:

It comes from the dimensions of virtual fundamental cycles.
After working out the corresponding formula I get that
one should define a new $\QQ$-grading on $HS(X)$.

If $Z$ is a connected component of $X_1$ we will define the rational number
$$age(Z)=\frac{1}{2 \pi i} \sum(\mbox{log eigenvalues of the action of $g$ on $Tx$}),$$
where $(x,g)$ is any point from $Z$.

The new $\QQ$-degree of the component  $H^k(Z)$ of $HS(X)$ I define as
$$k\ \ +\ \ 2*age.$$

Notice that the $\ZZ/2$ grading (superstructure)
is the old one, and has now nothing to do with the $Q$-grading.
Poincar\'e duality on $HS$ comes from dualities on each $H(Z)$.

That's all, and I am a bit tired of typing by now.

Please write me if something is not clear.

There are some useful examples:

$Y=2$dim torus, $G=Z/2$ with antipodal action.

Then $rk(HS)=6$.

Another case: $Y=($complex surface $S)^n$, $G=$symmetric group $S_n$.

$HS(Y)=H($Hilbert scheme of $S$\\ \hspace*{1in} resolving singularities
in the symmetric power of $S)$\footnote{Isomorphisms as  groups, the relationship of ring structures is delicate - D.A.}.

All the best,

Maxim Kontsevich
\subsubsection{Letter of  31 July 1996}\hfill

\noindent
Date: Wed, 31 Jul 96 11:21:35 +0200 \\
From: maxim@ihes.fr (Maxim Kontsevich)\\
To: lborisov@msri.org\\
Subject: Re: stable maps to orbifolds\\

Dear Lev,

First of all, there was a misprint in the definition of $X_1$ as you noticed\footnote{Borisov pointed out that an evaluation map to the non-rigidified inertia stack $\iX$ (here $X_1$), and even just to $\cX$, may not exist. Kontsevich defines it in what follows on the level of coarse moduli spaces, which works fine -  D.A.}.

I forgot in my letter to give the definition of the evaluation map.
Namely, if $S_1$ is an orbit of $G$ acting on the curve, we chose first a point
$p$ from $S_1$. The stabilizer $A$ of $p$ is a cyclic group because it acts
freely on a punctured neighborhood of  $p$. Also, this cyclic group has
a canonical generator g which rotates the tangent space to p to the minimal
angle in the anti-clockwise direction. I associate with the orbit $S_1$
the point $(\phi(p),g)$ modulo $G$ in $X_1$. This point is independent of the choice
of $p$ in $S_1$.

Best,

Maxim Kontsevich

P.S. Please call me Maxim, not Professor


\begin{thebibliography}{[20]}

\bibitem{Abramovich-Corti-Vistoli} D. 
 Abramovich, A. Corti and A. Vistoli, {\em  Twisted bundles and admissible covers.} Special issue in honor of Steven L. Kleiman.  Comm. Algebra  31  (2003),  no. 8, 3547--3618.

\bibitem{AGV0} D. Abramovich, T. Graber and A. Vistoli, {\em Algebraic orbifold quantum products.} Orbifolds in mathematics and physics (Madison, WI, 2001),  1--24, Contemp. Math., 310, Amer. Math. Soc., Providence, RI, 2002. 

\bibitem{AGV} --, {\em Gromov--Witten theory of Deligne--Mumford stacks}, preprint.

\bibitem{AVfibered} D. Abramovich and A. Vistoli, {\em  Complete moduli for fibered surfaces.}
Recent progress in intersection theory (Bologna, 1997), 1--31,
Trends Math.,
Birkh\"auser Boston, Boston, MA, 2000.
 
\bibitem{AV} --, {\em Compactifying the space of stable maps.}  J. Amer. Math. Soc.  15  (2002),  no. 1, 27--75. 


\bibitem{Behrend} K. Behrend, {\em 
Gromov-Witten invariants in algebraic geometry.} (English. English summary)
Invent. Math. 127 (1997), no. 3, 601--617.

\bibitem{Behrend-Fantechi} K. Behrend, and B. Fantechi, {\em  The intrinsic normal cone.}  Invent. Math.  128  (1997),  no. 1, 45--88.

\bibitem{Borisov-Chen-Smith} L. Borisov, L. Chen and G. Smith, {\em  The orbifold Chow ring of toric Deligne-Mumford stacks.}
J. Amer. Math. Soc. 18 (2005), no. 1, 193--215

\bibitem{BKR} T. Bridgeland,  A. King, M.  Reid, {\em Mukai implies McKay:
The McKay correspondence as an equivalence of derived categories.}
J. Amer. Math. Soc. 14 (2001), no. 3, 535--554 

\bibitem{Bryan-Graber-Pandharipande}, J. Bryan, T. Graber, and R. Pandharipande, {\em The orbifold quantum cohomology of $C^2/Z_3$ and Hurwitz-Hodge integrals}, preprint {\tt math.AG/0510335} 

\bibitem{Cadman} C. Cadman {\em Using stacks to impose tangency conditions on curves}, preprint {\tt math.AG/0312349}

\bibitem{Cadman2} C. Cadman {\em On the enumeration of rational plane curves with tangency conditions}, preprint {\tt math.AG/0509671}

\bibitem{COGP} P. 
Candelas, X.  de la Ossa, P. Green,  and L.  Parkes, {\em 
A pair of Calabi-Yau manifolds as an exactly soluble superconformal theory.}
Nuclear Phys. B 359 (1991), no. 1, 21--74.

\bibitem{Chen-Hu}  Chen B. and  Hu S., {\em A deRham model for Chen-Ruan cohomology ring of abelian orbifolds}, preprint {\tt math.SG/0408265}

\bibitem{Chen-Ruan} Chen W. and  Ruan Y., 
{\em Orbifold Gromov-Witten theory.}
Orbifolds in mathematics and physics (Madison, WI, 2001), 25--85,
Contemp. Math., 310,
Amer. Math. Soc., Providence, RI, 2002. 

\bibitem{DHVW} L. 
Dixon, J.A. Harvey, C.  Vafa, E.  Witten, {\em 
Strings on orbifolds.}
Nuclear Phys. B 261 (1985), no. 4, 678--686.


\bibitem{DHVW2} L. 
Dixon, J.A. Harvey, C.  Vafa, E.  Witten, {\em 
Strings on orbifolds II.}
 Nuclear Phys. B  274  (1986),  no. 2, 285--314.

\bibitem{Edidin} D. Edidin, {\em 
    Notes on the construction of the moduli space of curves}, preprint  {\tt math.AG/9805101}

\bibitem{Ekedahl} T. Ekedahl, {\em Boundary behaviour of Hurwitz schemes.}
The moduli space of curves (Texel Island, 1994), 173--198,
Progr. Math., 129,
Birkh\"auser Boston, Boston, MA, 1995.

\bibitem{Fantechi-Gottsche} B. Fantechi and L. G\"ottsche, { Orbifold cohomology for global quotients.}
Duke Math. J. 117 (2003), no. 2, 197--227.


\bibitem{Fulton-Pandharipande} W.\ Fulton  and R.\ Pandharipande,
{\em Notes on stable maps and quantum cohomology,} in {\it Algebraic
geometry---Santa Cruz 1995}, 45--96, Proc. Sympos. Pure Math., Part 2,
Amer. Math. Soc., Providence, RI, 1997.


\bibitem{Gillet} H.\ Gillet, {\em Intersection theory on algebraic stacks
and Q-varieties,} J. Pure Appl. Algebra 34, 193-240 (1984).

\bibitem{Givental} A. Givental, {\em 
Equivariant Gromov-Witten invariants.}
Internat. Math. Res. Notices 1996, no. 13, 613--663.

\bibitem{Goldin-Holm-Knutson} R. Goldin, T. Holm and A. Knutson, {\em Orbifold cohomology of torus quotients}, preprint {\tt math.SG/0502429}


\bibitem{Jarvis-Kaufmann-Kimura}
T. Jarvis, R. Kaufmann, and T. Kimura, {\em  Pointed admissible $G$-covers and $G$-equivariant cohomological field theories.}  Compos. Math.  141  (2005),  no. 4, 926--978.

\bibitem{JKK2}
T. Jarvis, R. Kaufmann, and T. Kimura, {\em Stringy K-theory and the Chern character}, preprint {\tt math.AG/0502280}

\bibitem{Keel-Mori} S. Keel and S. Mori, 
{\em Quotients by groupoids,} Ann. of
Math. (2) {\bf 145} (1997), no.~1, 193--213. 

\bibitem{Kontsevich} M.  Kontsevich, {\em Enumeration of rational curves via torus actions.}  The moduli space of curves (Texel Island, 1994),  335--368, Progr. Math., 129, Birkh\"auser Boston, Boston, MA, 1995.

\bibitem{Kontsevich-motivic} M. Kontsevich, {\em Grothendieck ring of motives and related rings}, notes for the lecture ``String Cohomology" at Orsay, December 7, 1995.\\ {\tt http://www.mabli.org/jet-preprints/Kontsevich-MotIntNotes.pdf}

\bibitem{Kontsevich-Manin}
M. Kontsevich and Yu. Manin, {\em 
Gromov-Witten classes, quantum cohomology, and enumerative geometry. }
Comm. Math. Phys. 164 (1994), no. 3, 525--562.

\bibitem{Li-Qin-Wang} Li W.-P., Qin Z., and Wang W.,    {\em  The cohomology rings of Hilbert schemes via Jack polynomials,} preprint {\tt math.AG/0411255}

\bibitem{Lian-Liu-Yau}
 Lian B.,  Liu K.,    Yau S.-T., {\em 
Mirror principle I}
Asian J. Math. 1 (1997), no. 4, 729--763.

\bibitem{Mochizuki} S. Mochizuki, {\em 
Extending families of curves over log regular schemes.} J. Reine Angew. Math. 511 (1999), 43--71.

\bibitem{Olsson-stacks} M. Olsson, {\em Sheaves on Artin stacks}, preprint \\ {\tt http://www.ma.utexas.edu/users/molsson/qcohrevised.pdf}

\bibitem{Olsson-Hom} --, {\em Hom--stacks and restriction of scalars}, preprint \\ {\tt http://www.ma.utexas.edu/users/molsson/homstackfinal.pdf}


\bibitem{Olsson-log} --, {\em On (log) twisted curves}, preprint \\ {\tt http://www.ma.utexas.edu/users/molsson/Logcurves.pdf}

\bibitem{Romagny} M. Romagny, {\em 
Group actions on stacks and applications.}
Michigan Math. J. 53 (2005), no. 1, 209--236.

\bibitem{Tseng}  Tseng H.-H., {\em Orbifold Quantum Riemann-Roch, Lefschetz and Serre}, preprint {\tt math.AG/0506111}

\bibitem{Uribe}  B. Uribe, {\em 
Orbifold cohomology of the symmetric product. }
Comm. Anal. Geom. 13 (2005), no. 1, 113--128.


\bibitem{Vistoli} A. Vistoli, {\em Intersection theory on
algebraic 
stacks and on their moduli spaces.} Invent. Math. 97 (1989), no. 3,
613-670.

\bibitem{Zaslow} E. Zaslow, {\em Topological orbifold models and quantum cohomology rings.} Comm. Math. Phys. 156 (1993), no. 2, 301--331.

\end{thebibliography}
\end{document}